\documentclass[final]{siamltex}

\usepackage{amsfonts}
\usepackage{graphicx}
\usepackage{amsmath}
\usepackage{graphicx}
\usepackage{amssymb}
\usepackage[colorlinks,linkcolor={blue}]{hyperref}%
\newtheorem{remark}[theorem]{Remark}

\newcommand {\integers}{{\mathbb{Z}}}   
\newcommand {\reals}{{\mathbb{R}}}      
\newcommand {\complexs}{{\mathbb{C}}}   
\newcommand {\tori}{{\mathbb{T}}}       

\usepackage{color}    
\definecolor{Red}{rgb}{1,0,0}

\usepackage{ulem}    
\usepackage{url}    

\DeclareGraphicsExtensions{.eps,.pdf,{}}

\title{Existence of a Center Manifold in a Practical Domain around $L_1$ in
the Restricted Three Body Problem}

\author{Maciej J. Capi\'nski ({\tt mcapinsk@agh.edu.pl})\thanks{Supported by the Polish State Ministry of Science and Information Technology grant N201 543238. The work was initiated during a visit of the first author to the 2008 SIMS workshop in Barcelona, sponsored by the Centre de Recerca Matem\`atica. 
The work was concluded during a visit of the first author to the University of Texas at Austin sponsored by the Ko\'sciuszko Foundation.} \and 
Pablo Rold\'an ({\tt pablo.roldan@upc.edu})\thanks{Pablo Rold\'an was
supported in part by MICINN-FEDER grant MTM2009-06973 and CUR-DIUE grant
2009SGR859.}}

\begin{document}
	\maketitle
\begin{abstract}
We give a proof of existence of centre manifolds within large domains for systems with an integral of motion. The proof is based on a combination of topological tools, normal forms and rigorous-computer-assisted computations. We apply our method to obtain an explicit region in which we prove existence of a center manifold in the planar Restricted Three Body Problem.  
\end{abstract}

\begin{keywords}
center manifolds, normal forms, 
celestial mechanics, restricted
three-body problem, covering relations, cone conditions.
\end{keywords}

\begin{AMS}
37D10,    	
37G05,    	
37N05,    	
34C20,    	
34C45,    	
70F07,    	
70F15,   	
70K45.    	
\end{AMS}

\pagestyle{myheadings}
\thispagestyle{plain}
\markboth{M. J. CAPI\'NSKI AND P. ROLD\'AN}{CENTER MANIFOLDS IN A PRACTICAL
DOMAIN}



\section{Introduction}

Center manifolds are an important tool for the local analysis of dynamical
systems. In this paper we develop a methodology to prove the existence of
center manifolds in a ``large'' neighbourhood of the equilibrium point. The method involves the use of normal forms, topological results, and computer assisted
computations. The novelty of our approach is that it provides \emph{explicit
rigorous bounds} on the size and location of the manifold for a given
dynamical system. Moreover, under appropriate hypothesis we prove that the
manifold is unique.

In contrast, the classical center manifold theorems show existence of a
manifold in \emph{some} neighbourhood, but they do not readily provide
information on the size of this neighbourhood. Also, the classical normal form
theorems construct an accurate approximation to the dynamics in a
neighbourhood, but the normal form is usually not convergent. Sometimes the
normal form does converge, but we lack information on its domain of convergence.

To show the power of our methodology, in the second part of this paper we
prove existence and uniqueness of the center manifold in a practical domain
around an equilibrium point of the celebrated Restricted Three Body Problem
(RTBP). By practical we mean that such domain possibly could be used for
realistic space mission design, since it is not too small. To our knowledge,
this is the first proof of existence of the center manifold in a practical
domain for the RTBP.

For the rest of this introduction, we define the center manifold and mention
some previous results related to this paper. Finally we explain how this paper
is organized into sections.

\begin{definition}
Consider a differential equation on $\reals^{n}$
\begin{equation}
\dot{x}=Ax+f(x),\label{eq:diffeq}%
\end{equation}
where $A$ is linear and $f$ has no constant or linear terms. The origin is a
fixed point. Let $\reals^{n}=E^{c}\oplus E^{u}\oplus E^{s}$ be the usual
decomposition into the center, unstable, and stable invariant subspaces with
respect to $A$.

A \emph{center manifold} $W^{c}$ is an invariant manifold of the flow of
\eqref{eq:diffeq}, tangent to $E^{c}$ at the origin, and of the form
\[
W^{c}=\{(\theta,\chi(\theta))\colon\ \theta\in U\},
\]
where $\chi\colon\ E^{c}\rightarrow E^{u}\oplus E^{s}$ is a $C^{k}$ function,
and $U$ is an open neighborhood of $0$ in $E^{c}$.
\end{definition}

We are naturally lead to study the flow in the center manifold. The center
manifold approach has the advantage that this reduced problem is a dynamical
system on a lower-dimensional manifold (of the same dimension as $E^{c}$). The
reduced problem contains crucial information of the full
problem~\eqref{eq:diffeq}. The qualitative behavior of the flow on the center
manifold completely determines the behavior of the full flow around the fixed point \cite{Car}. Also,
every center manifold contains all globally bounded solutions (e.g. fixed
points and periodic orbits) which are close enough to the
origin \cite{Sijbrand}.

Let us now mention some results related to this paper. The existence of center
manifolds is discussed in many dynamical systems books, for instance in
\cite{GH, ChowHale}, and the monograph~\cite{Car}. The subtle
properties of center manifolds such as (non)-uniqueness, tangency, limited
differentiability, and (non)-analyticity are discussed in~\cite{Sijbrand}.

Lyapunov~\cite{Lyapunov} studied the case in which the linear operator $A$ in
equation~\eqref{eq:diffeq} has a simple pair of eigenvalues $\pm\omega i$. He
proved the existence of a center manifold filled with an analytic
one-parameter family of periodic orbits. The main hypothesis are the presence
of an integral of motion and a nonresonance condition. Such situation arises
for the equilibrium point of the Restricted Three Body Problem that we study
in the second part of this paper. Lyapunov's theorem applies to the
Restricted Three Body Problem (cf.~\cite{SM} \S 18), but again is only local
and does not readily provide estimates on its domain of validity.

Normal forms make a very powerful and general technique to approximate local
dynamics, including the center manifold, and stable/unstable manifolds of a
fixed point. It is also a classical subject discussed in many dynamical system
books, for instance \cite{MeyerHall} (Hamiltonian systems), \cite{ChowHale} (general
differential equations), \cite{GH}, and the monograph~\cite{Murdock}.

Given their usefulness, normal forms have been applied to approximate the
center manifold around the equilibrium points of the planar Restricted Three
Body Problem \cite{CanaliasMasdemont,CDMR}, and the spatial RTBP
\cite{JorbaM99, DelshamsMR08}. In particular, our implementation of normal forms is
based on~\cite{Jorba}. This technique has important applications in space
mission design \cite{GJSM, GKLMMR} and diffusion estimates
\cite{JorbaSimo94, JorbaVillanueva98}.

Regarding the planar RTBP, we would also like to mention the numerical
explorations of Broucke \cite{Broucke}, where he performed an extensive study of
different families of periodic orbits. In particular, he finds a family of
numerical periodic orbits around the same equilibrium point that we study in
this paper. The family extends up to a very large neighborhood of the
equilibrium point (much larger than our rigorous result).

The paper is organized as follows. In Section \ref{sec:setup} we give the
setup of the problem and state our main theorem (Theorem \ref{th:main}). 
Assumptions of the theorem are based on estimates on the derivatives of the vector field 
within the investigated region. Based on these the existence of an invariant manifold is established. In
Section \ref{sec:top-maps} we give a topological proof of the existence of an
invariant manifold for maps with saddle-center-type properties. In Section
\ref{sec:proof-of-main} we use the result obtained for maps to prove 
Theorem \ref{th:main}. In Section \ref{sec:3bp} we apply our Theorem \ref{th:main} to prove the
existence of a center manifold around an equilibrium point $L_{1}$ in the
RTBP. To do so we first introduce the problem and present a procedure of transforming the 
system into a normal form. We then discuss how normal forms
provide very accurate approximations of center manifolds. Finally we combine Theorem \ref{th:main} and normal forms with
rigorous interval arithmetic based computer assisted computations to prove the
existence of the manifold. Section \ref{sec:concluding-remarks} contains
concluding remarks and an outline of future work.



\section{Setup}

\label{sec:setup}

We will consider the following problem. Let $F:\mathbb{R}^{n}\rightarrow
\mathbb{R}^{n}$ and
\begin{equation}
\mathbf{x}^{\prime}=F(\mathbf{x}) \label{eq:our-ode}%
\end{equation}
be an ODE (we impose the usual assumptions implying existence and uniqueness
of solutions) with a fixed point $\mathbf{x}_{0}$ and an integral of motion
$H:\mathbb{R}^{n}\rightarrow\mathbb{R}$. By this we mean that for any solution
$q(t)$ of (\ref{eq:our-ode}) we have%
\begin{equation}
H(q(t))=c, \label{eq:energy-cond}%
\end{equation}
where $c$ is some constants dependent on the initial condition $q(0)$. Since
in our applications we shall deal with the restricted three body problem,
which is a Hamiltonian system where $H$ is the Hamiltonian, we
shall refer to $H$ as the energy from now on. We shall use a notation
$\Phi(t,\mathbf{x})$ for the flow induced by (\ref{eq:our-ode}).

\subsection{Well aligned coordinates \label{sec:adjusted-coordinates}}

We will investigate the dynamics of (\ref{eq:our-ode}) in some compact set
$D$, contained in an open subset $U$ of $\reals ^n$, such that the fixed point $\mathbf{x}_0\in D$, and whose
image by a diffeomorphism
\begin{equation}
\phi:U\rightarrow\phi(U)\subset\mathbb{R}^{n} \label{eq:phi-total-change}%
\end{equation}
is%
\begin{equation}
\phi(D)=D_{\phi}=\bar{B}_{c}^{R}\times\bar{B}_{u}^{r}\times\bar{B}_{s}^{r},
\label{eq:Dphi-r}%
\end{equation}
where $\bar{B}_{i}^{r}$ (for
$i\in\{c,u,s\}$) stand for $i$-dimensional closed balls around zero of radius
$r$. We assume that $n=c+u+s$. We will refer to $p=\phi(\mathbf{x})$ as the
\emph{aligned coordinates}. In these coordinates we will use a notation
$p=(\theta,x,y)$ with $\theta\in\bar{B}_{c}^{R},$ $x\in\bar{B}_{u}^{r}$ and
$y\in\bar{B}_{s}^{r}$. We will refer to $\theta$ as the central coordinate, to
$x$ as the unstable coordinate and to $y$ as the stable coordinate (the
subscripts $c,u,s$ standing for central, unstable and stable respectively).

The motivation behind the above setup is the following. We will search for a
center manifold of (\ref{eq:our-ode}) homeomorphic to a $c$-dimensional disc
inside the set $D$. Such manifolds have associated stable and unstable vector
bundles, which in the coordinate system $\phi$ are given approximately by the
coordinates of the balls $\bar{B}_{s}^{r}$ and $\bar{B}_{u}^{r}$ respectively.
We do not assume though that the coordinates $x$ and $y$ align exactly with
directions of hyperbolic expansion and contraction. It will turn out that it
is enough that they point roughly in these directions. The remaining coordinates
$\theta$ are the central coordinates of our system. We need to have a good
guess on where the center manifold is. This guess is given by $\phi^{-1}(\bar{B}_{c}^{R}%
\times\{0\})\subset\mathbb{R}^{n}$.  The change of coordinates $\phi$ can be obtained from some non-rigorous numerical computation (in our application for the RTBP - normal forms). 
It is important to emphasise that
we will not assume that $\phi^{-1}(\bar{B}_{c}^{R}\times\{0\})$ is invariant
under the flow (\ref{eq:our-ode}). Allowing for errors, we expect the true manifold to lie in $\phi^{-1}(\bar
{B}_{c}^{R}\times\bar{B}_{u}^{r}\times\bar{B}_{s}^{r}).$ This means that we
take an enclosure of radius $r$ of our initial guess and look for the
invariant manifold in this enclosure.

We will search for the part of the center manifold with energy $H\leq h$ for
some $h\in\mathbb{R}$. We assume that the center coordinate is well aligned
with the energy $H$ in the sense that we have%
\begin{equation}
H(\phi^{-1}(\bar{B}_{c}^{R-v}\times\bar{B}_{u}^{r}\times\bar{B}_{s}%
^{r}))<h<H\left(  \phi^{-1}(\partial\bar{B}_{c}^{R}\times\bar{B}_{u}^{r}%
\times\bar{B}_{s}^{r})\right)  , \label{eq:h-boundary-bound}%
\end{equation}
for some $v>0$ (here we use a notation $\partial A$ to denote the boundary of
a set $A$).

Our detection of the center manifold in the RTBP is going to be carried out in two stages.
First we obtain $\phi$ as a change of coordinates into a normal form, after which we shall employ our topological theorem (Theorem \ref{th:main}) to prove the existence of the manifold.

\subsection{Local bounds on the vector field and the statement of the main
result}

\label{sec:main-statement}

We are now ready to state the main assumptions needed for our method. These
will be expressed in terms of local bounds on the derivative of the vector
field (\ref{eq:our-ode}). First let us introduce a notation $F^{\phi}$ for the
vector field in the aligned coordinates i.e.%
\begin{equation}
F^{\phi}(p)=D\phi(\phi^{-1}(p))F(\phi^{-1}(p)), \label{eq:field-local}%
\end{equation}
and a notation $[dF^{\phi}(N)]$ for an interval enclosure of the derivative on
a set $N\subset D_{\phi}$%
\[
\lbrack dF^{\phi}(N)]=\left\{  A\in\mathbb{R}^{n\times n}|A_{ij}\in\left[
\inf_{p\in N}\frac{dF_{i}^{\phi}}{dp_{j}},\sup_{p\in N}\frac{dF_{i}^{\phi}%
}{dp_{j}}\right]  ,\text{ for all }i,j=1,\ldots,n\right\}  .
\]

For any point $p=(\theta,0,0)$ from $\bar{B}_{c}^{R}\times\{0\}\times\{0\}$ we
define a set%
\begin{equation}
N_{p}:=\bar{B}_{c}(\theta,\rho)\times\bar{B}_{u}^{r}\times\bar{B}_{s}^{r}\cap
D_{\phi},\label{eq:Nq-set}%
\end{equation}
where $\bar{B}_{c}(\theta,\rho)$ is a $c$ dimensional ball of radius
$\rho>0\ $centred at $\theta$. We introduce the following notations for the
bound on the derivatives of $F^{\phi}$ on the sets $N_{p}$%
\begin{equation}
\lbrack dF^{\phi}(N_{p})]\subset\left(
\begin{array}
[c]{ccc}%
\mathbf{C} & \boldsymbol{\varepsilon}_{c} & \boldsymbol{\varepsilon}_{c}\\
\boldsymbol{\varepsilon}_{m} & \mathbf{A} & \boldsymbol{\varepsilon}_{u}\\
\boldsymbol{\varepsilon}_{m} & \boldsymbol{\varepsilon}_{s} & \mathbf{B}%
\end{array}
\right)  .\label{eq:dF-bound}%
\end{equation}
Here $\mathbf{A,B,C,}\boldsymbol{\varepsilon}_{c},\boldsymbol{\varepsilon}%
_{m},\boldsymbol{\varepsilon}_{s}$ and $\boldsymbol{\varepsilon}_{u}$ are
interval matrices, that is matrices with interval coefficients. Here we
slightly abuse notations since the pairs of matrices $\boldsymbol{\varepsilon
}_{c}$ and $\boldsymbol{\varepsilon}_{m}$ need not be equal; they even have
different dimension when $u\neq s$. We use the same notation since later on we
shall assume uniform bounds for both of matrices $\boldsymbol{\varepsilon}%
_{c}$ and both $\boldsymbol{\varepsilon}_{m}$. Let us also note that the
bounds $\mathbf{A,B,C,}\boldsymbol{\varepsilon}_{c},\boldsymbol{\varepsilon
}_{m},\boldsymbol{\varepsilon}_{s}$ and $\boldsymbol{\varepsilon}_{u}$ may be
different for different $p.$ We do not indicate this in our notations to keep
them relatively simple.

\begin{remark}
If the system possesses a center manifold and the adjusted coordinates are
well aligned in the sense of section \ref{sec:adjusted-coordinates}, then the
interval matrices $\boldsymbol{\varepsilon}_{i}$ in (\ref{eq:dF-bound}), with
$i\in\{c,m,s,u\}$ should turn out to be small. The matrices $\mathbf{A,B,C}$
are the bounds on derivatives of the vector field in the unstable, stable and
central directions respectively. If the alignment of our coordinates is
correct then we expect the contraction/expansion rates associated with
$\mathbf{C}$ to be weaker than for $\mathbf{A}$ and $\mathbf{B.}$
\end{remark}

We will use the following notations to express our assumptions about
$[dF^{\phi}(N_{p})].$ Let $\delta^{u},\delta^{s},c^{u},c^{s},\varepsilon_{i} >
0$ denote contraction/expansion rates, such that for any matrix $A\in
\mathbf{A}$, $B\in\mathbf{B,}$ $e_{i}\in\boldsymbol{\varepsilon}_{i}$ for
$i\in\{m,c,u,s\}$, we have%

\begin{gather}
\inf\{x^{T}Ax:||x||=1\}>\delta^{u},\label{eq:est-expansion}\\
\sup\{y^{T}By:||y||=1\}<-\delta^{s},\label{eq:est-contraction}\\
c^{s}<\inf\{\theta^{T}C\theta:||\theta||=1\}\leq\sup\{\theta^{T}%
C\theta:||\theta||=1\}<c^{u},\label{eq:est-central}\\
||e_{i}||<\varepsilon_{i}\quad\text{for }i\in
\{m,c,u,s\}.\label{eq:est-epsilons}%
\end{gather}
Once again, $\varepsilon_{i},c^{s},c^{u},\mu,\delta^{u}$ and $\delta^{s}$ can
depend on $p$.

Let $\gamma,\alpha_{h},\alpha_{v},\beta_{h},\beta_{v}>0$ be constants such
that
\begin{equation}
\alpha_{h}>\alpha_{v}\quad\text{and}\quad\beta_{v}>\beta_{h}%
,\label{eq:alpha-beta-setup}%
\end{equation}
and such that the radius $\rho$ considered for the central part of the sets
$N_{p}$ satisfies%
\begin{equation}
\rho>r\sqrt{\frac{\alpha_{h}}{\gamma}},\quad\quad\rho>r\sqrt{\frac{\beta_{v}%
}{\gamma}},\qquad\label{eq:rho-cond-1}%
\end{equation}
where $r$ is the radius of the balls $\bar{B}_{u}^{r}$ and $\bar{B}_{s}^{r}$
in (\ref{eq:Dphi-r}). Let us define the following constants%
\begin{align}
\kappa_{c}^{\text{forw}} &  :=c^{u}+\frac{1}{2}\left(  \frac{\alpha_{h}%
}{\gamma}\varepsilon_{m}+\frac{\beta_{h}}{\gamma}\varepsilon_{m}%
+2\varepsilon_{c}\right)  ,\nonumber\\
\kappa_{u}^{\text{forw}} &  :=\delta^{u}-\frac{1}{2}\left(  \varepsilon
_{m}+\varepsilon_{u}+\frac{\gamma}{\alpha_{h}}\varepsilon_{c}+\frac{\beta_{h}%
}{\alpha_{h}}\varepsilon_{s}\right)  ,\label{eq:kappa-h}\\
\kappa_{s}^{\text{forw}} &  :=-\delta^{s}+\frac{1}{2}\left(  \varepsilon
_{m}+\frac{\alpha_{h}}{\beta_{h}}\varepsilon_{u}+\frac{\gamma}{\beta_{h}%
}\varepsilon_{c}+\varepsilon_{s}\right)  ,\nonumber
\end{align}%
\begin{align}
\kappa_{c}^{\text{back}} &  :=c^{s}-\frac{1}{2}\left(  \varepsilon_{m}%
\frac{\alpha_{v}}{\gamma}+\varepsilon_{m}\frac{\beta_{v}}{\gamma}%
+2\varepsilon_{c}\right)  ,\nonumber\\
\kappa_{u}^{\text{back}} &  :=\delta^{u}-\frac{1}{2}\left(  \varepsilon
_{m}+\varepsilon_{u}+\frac{\gamma}{\alpha_{v}}\varepsilon_{c}+\frac{\beta_{v}%
}{\alpha_{v}}\varepsilon_{s}\right)  ,\label{eq:kappa-v}\\
\kappa_{s}^{\text{back}} &  :=-\delta^{s}+\frac{1}{2}\left(  \varepsilon
_{m}+\frac{\alpha_{v}}{\beta_{v}}\varepsilon_{u}+\frac{\gamma}{\beta_{v}%
}\varepsilon_{c}+\varepsilon_{s}\right)  .\nonumber
\end{align}
The superscripts "forw" and "back" in the above constants come from the fact
that they shall be associated with estimates on the dynamics induced by the
vector field (\ref{eq:field-local}), for forward and backward evolution in
time respectively. At this stage the subscripts $v$ and $h$ in constants
$\alpha$ and $\beta$ do not have an intuitive meaning. During the course of
the proof they shall be associated with horizontal and vertical slopes of
constructed invariant manifolds (hence $h$ for "horizontal" and $v$ for
"vertical"), and then their meaning will become more natural.

\begin{remark}
Even though coefficients (\ref{eq:kappa-h}), (\ref{eq:kappa-v}) are technical in nature, they have a quite
natural interpretation in terms of the dynamics of the system. The estimates
$\kappa_{c}^{i},\kappa_{u}^{i},\kappa_{s}^{i}$ for $i\in\{$forw,back$\}$ are
essentially estimates on the contraction/expansion rates associated with the
center, unstable and stable coordinates respectively. These estimates take
into account errors $\varepsilon_{i}$ for $i\in\{s,u,c,m\}$ in the setup of
coordinates. Note that when our coordinates are perfectly aligned with the
dynamics, then $\varepsilon_{i}=0$ for $i\in\{s,u,c,m\},$ and in turn%
\[
\kappa_{s}^{\text{forw}}=\kappa_{s}^{\text{back}}=-\delta^{s},\quad\quad
\kappa_{c}^{\text{back}}=c^{s},\quad\quad\kappa_{c}^{\text{forw}}=c^{u}%
,\quad\quad\kappa_{u}^{\text{forw}}=\kappa_{u}^{\text{back}}=\delta^{u},
\]
which are the bounds on the derivative of the vector field in the unstable,
stable and center directions given in (\ref{eq:est-expansion}),
(\ref{eq:est-contraction}), (\ref{eq:est-central}). The key assumptions of
Theorem \ref{th:main} are (\ref{eq:mth-as-cone1}) and (\ref{eq:mth-as-cone2}).
In particular, these assumptions imply
\[
\kappa_{s}^{\text{back}}<\kappa_{c}^{\text{back}}\quad\quad\kappa
_{c}^{\text{forw}}<\kappa_{u}^{\text{forw}},
\]
which is equivalent to assuming that the dynamics in the center
coordinate is weaker than dynamics in the stable and unstable directions.
These are classical assumptions for center manifold theorems (See \cite{GH},
for instance).
\end{remark}

\begin{remark}
We have certain freedom of choice for the constants $\gamma$, $\alpha_{h}$%
, $\alpha_{v}$, $\beta_{h}$, $\beta_{v}$. They offer flexibility when verifying
assumptions of Theorem \ref{th:main}. During the course of the proof of
Theorem \ref{th:main} it will turn out that they also give Lipschitz type
bounds $L_{c}=\sqrt{\frac{2\gamma}{\min(\alpha_{h}-\alpha_{v},\beta_{v}%
-\beta_{h})}},$ $L_{s}=\sqrt{\frac{1}{\alpha_{h}}\max(\gamma,\beta_{h})},$
$L_{u}=\sqrt{\frac{1}{\beta_{v}}\max(\gamma,\alpha_{v})}$ for our center,
stable and unstable manifolds respectively (for more details see Corollary
\ref{cor:lip-bounds}).
\end{remark}

We are now ready to state our main tool for detection of center manifolds.

\begin{theorem}
\label{th:main}(Main Theorem) Let $h\in\mathbb{R}.$ Assume that
(\ref{eq:h-boundary-bound}) holds for some $v>0.$ Assume also that for any
$p\in\bar{B}_{c}^{R}\times\{0\}\times\{0\},$ for the constants $\kappa
_{c}^{\text{forw}}$, $\kappa_{u}^{\text{forw}}$, $\kappa_{s}^{\text{forw}}$%
, $\kappa_{c}^{\text{back}}$, $\kappa_{u}^{\text{back}}$, $\kappa_{s}^{\text{back}%
}$, $\varepsilon_{u}$, $\varepsilon_{s}$, $\delta^{u}$, $\delta^{s}$ computed on a set
$N_{p}$ (defined by (\ref{eq:Nq-set})) the following inequalities hold:%
\begin{align}
\kappa_{c}^{\text{forw}},\kappa_{s}^{\text{forw}} &  <\kappa_{u}^{\text{forw}%
},\hspace{2cm} 0<\kappa_{u}^{\text{forw}},\label{eq:mth-as-cone1}\\
\kappa_{s}^{\text{back}} &  <0,\hspace{2cm}\kappa_{s}^{\text{back}%
}<\kappa_{c}^{\text{back}},\kappa_{u}^{\text{back}},\label{eq:mth-as-cone2}%
\end{align}
and also that there exist $E_{u},E_{s}>0$ such that for any $q\in N_{p}%
\cap(\bar{B}_{c}^{R}\times\{0\}\times\{0\})$
\begin{equation}
||\pi_{x}F^{\phi}(q)||<rE_{u},\quad\text{\quad}||\pi_{y}F^{\phi}%
(q)||<rE_{s},\label{eq:mth-as-cover-3}%
\end{equation}
and
\begin{align}
E_{u}+\varepsilon_{u} &  <\delta^{u},\label{eq:mth-as-cover-1}\\
E_{s}+\varepsilon_{s} &  <\delta^{s}.\label{eq:mth-as-cover-2}%
\end{align}
If above assumptions hold, then there exists a $C^{0}$ function%
\[
\chi:\bar{B}_{c}^{R-v}\rightarrow D_{\phi}%
\]
such that

\begin{enumerate}
\item For any $\theta\in\bar{B}_{c}^{R-v}$ we have $\pi_{\theta}\chi
(\theta)=\theta$ and
\[
\Phi(t,\phi^{-1}(\chi(\theta)))\in D\quad\text{for all }t\in\mathbb{R}.
\]

\item If for some $\mathbf{x}\in\phi^{-1}(\bar{B}_{c}^{R-v}\times\bar{B}%
_{u}^{r}\times\bar{B}_{s}^{r})$ we have
\[
\Phi(t,\mathbf{x})\in D\quad\text{for all }t\in\mathbb{R}%
\]
then there exists a $\theta\in\bar{B}_{c}^{R-v}$ such that $\mathbf{x}%
=\phi^{-1}(\chi(\theta))$.
\end{enumerate}
\end{theorem}

In subsequent sections we shall present a proof of this theorem building up
auxiliary results along the way. Before we move on to these results let us
make a couple of comments on the result.

\begin{remark}
Theorem \ref{th:main} establishes uniqueness of the invariant manifold. This
is not a typical scenario in case of center manifolds which are usually not
 unique. Uniqueness in our case follows from condition
(\ref{eq:h-boundary-bound}), which by our construction will ensure that for
any point from our center manifold a trajectory starting from it cannot leave
the set $D$. This means that dynamics on the center manifold with $H\leq h$ is
contained in a compact set. This is the underlying reason that allows us to
obtain uniqueness.
\end{remark}

\begin{remark}
The main strength of our result lies in the fact that it allows us to easily
obtain explicit bounds for the position and size of the manifold. The center
manifold is contained in $D=\phi^{-1}(\bar{B}_{c}^{R}\times\bar{B}_{u}%
^{r}\times\bar{B}_{s}^{r})$. Since the manifold is a graph of $\chi,$ from
point 1. of Theorem \ref{th:main} we know that it is of the form $\phi
^{-1}\left\{  \left(  \theta,\pi_{x,y}\chi(\theta)\right)  |\theta\in\bar
{B}_{c}^{R-v}\right\}  ,$ which ensures that it "fills in" the set $D$ nontrivially.
In contrast, the classical center manifold theorem does not provide such explicit
bounds. 
\end{remark}

\begin{remark}
In principle, one could derive some explicit analytic bounds using e.g. the
``method of majorands'' explained in the book of Siegel--Moser~\cite{SM}.
However, to guarantee existence of the center manifold in a neighborhood of
the equilibrium point that is not too small, one would require a substantial
amount of very careful estimates.
\end{remark}

\begin{remark}
It is important to remark that our result only establishes continuity
(together with Lipschitz type conditions) of the center manifold. The center
manifold theorem clearly indicates that in a sufficiently small neighbourhoods
of a saddle-center fixed point we should have higher order smoothness. We
believe though that similar in spirit assumptions to those of Theorem
\ref{th:main} should imply higher order smoothness. This will be the subject
of forthcoming work. The result obtained so far should be regarded as a first
step towards this end.

In our application for the RTBP, in a neighbourhood sufficiently close to the equilibrium point, our manifold shall inherit all regularity which follows from the center manifold theorem (see Remark \ref{rem:LapMoser}).
\end{remark}

Let us finish the section with a final comment. In order to verify assumptions
of Theorem \ref{th:main} it is sufficient to consider some finite covering
$\{\bigcup_{i\in I}U_{i}\}$ of the set $D_{\phi}$ and to verify bounds on
local derivatives on sets $U_{i}.$ It is not necessary to consider an infinite
number of points $p$ and their associated sets $N_{p},$ as long as for any
$p\in\bar{B}_{c}^{R}\times\{0\}\times\{0\}$ we have $N_{p}\subset U_{i}$ for
some $i\in I.$ This makes assumptions of Theorem \ref{th:main} verifiable in
practice using rigorous computer assisted tools.



\section{Topological approach to centre manifolds for maps}

\label{sec:top-maps}

In this section we will state some preliminary results, which will next be
used for the proof of Theorem \ref{th:main} in Section \ref{sec:proof-of-main}%
. The results will be stated for maps instead of flows. In Section
\ref{sec:proof-of-main} we will take a time shift along a trajectory map for
the flow generated by (\ref{eq:our-ode}) and apply the results to it. The main
result of this section is Theorem \ref{th:nhim-th-flows}. The result is in the
spirit of versions of normally hyperbolic invariant manifold theorems obtained
in \cite{Ca}, \cite{CZ} and \cite{CS}. The main difference is that we do not
deal with a normally hyperbolic manifold without boundary, but with a selected
part of a centre manifold (homeomorphic to a disc) with a boundary. In this
section the fact that the dynamics does not diffuse through the boundary along the
centre coordinate is imposed by assumption. This assumption will later follow
from assuming that (\ref{eq:energy-cond}), (\ref{eq:h-boundary-bound}) hold for (\ref{eq:our-ode}).

We now give the setup for maps. Let $D\subset U\subset\mathbb{R}^{n}$, the
change of coordinates $\phi:U\rightarrow\phi(U)$, and $D_{\phi}=\phi(D)$, be
as in Section \ref{sec:adjusted-coordinates}. We consider a dynamical system
given by a smooth invertible map $f:U\rightarrow U.$ In adjusted coordinates
we denote the map as $f_{\phi}:=\phi\circ f\circ\phi^{-1}$, $f_{\phi}%
:\phi(U)\rightarrow\mathbb{R}^{n}.$ We assume that
\begin{equation}
H(p)=H(f(p)) \label{eq:map-integral}%
\end{equation}
for all $p\in D_{\phi}$ and also that for some $v>0$ condition
(\ref{eq:h-boundary-bound}) holds.

We introduce the following sets
\begin{align}
D_{\phi}^{-} &  =\bar{B}_{c}^{R}\times\partial\bar{B}_{u}^{r}\times\bar{B}%
_{s}^{r},\label{eq:Dphi-minus}\\
D_{\phi}^{+} &  =\bar{B}_{c}^{R}\times\bar{B}_{u}^{r}\times\partial\bar{B}%
_{s}^{r}.\nonumber
\end{align}
We now introduce a number of definitions. The first is a definition of a
covering relation. \begin{figure}[ptb]
\begin{center}
\includegraphics[width=3in]{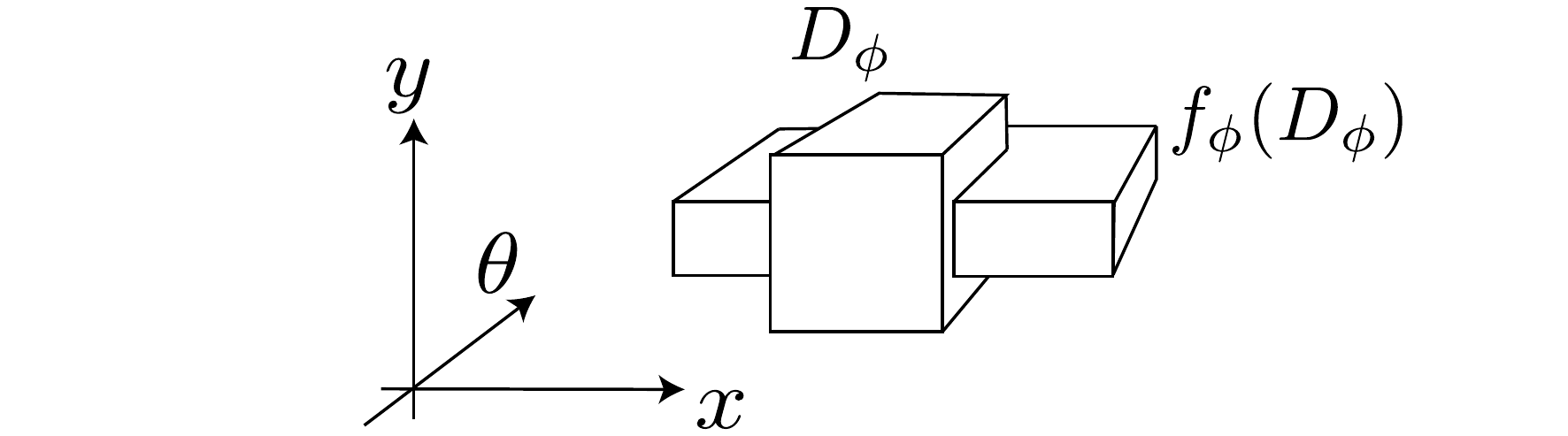}
\end{center}
\caption{A map $f$ which satisfies covering conditions. The set $D_{\phi}$ is
contracted in  coordinate $y$ and expanded in  coordinate $x$. Note that
in the $\theta$ coordinate the set may be simply shifted, expanded or
contracted, just as long as conditions (\ref{eq:cover-cond-1}%
),...,(\ref{eq:cover-cond-5}) are satisfied. }%
\label{fig:covering}%
\end{figure}

\begin{definition}
We say that a map $f:U\rightarrow U$ satisfies \emph{covering conditions} in
$D$ if
\begin{align}
\pi_{x}(f_{\phi}(D_{\phi}^{-}))\cap\bar{B}_{u}^{r} &  =\emptyset
,\label{eq:cover-cond-1}\\
\pi_{y}(f_{\phi}^{-1}(D_{\phi}^{+}))\cap\bar{B}_{s}^{r} &  =\emptyset
,\label{eq:cover-cond-2}\\
\pi_{y}(f_{\phi}(D_{\phi}))\cap\left(  \mathbb{R}^{s}\setminus\bar{B}_{s}%
^{r}\right)   &  =\emptyset,\label{eq:cover-cond-3}\\
\pi_{x}(f_{\phi}^{-1}(D_{\phi}))\cap\left(  \mathbb{R}^{u}\setminus\bar{B}%
_{u}^{r}\right)   &  =\emptyset,\label{eq:cover-cond-4}%
\end{align}
and for any point $p\in\bar{B}_{c}^{R}\times\{0\}$,%
\begin{equation}
\pi_{(x,y)}f_{\phi}(p),\,\pi_{(x,y)}f_{\phi}^{-1}(p)\in\mathrm{int}\left(
\bar{B}_{u}^{r}\times\bar{B}_{s}^{r}\right)  .\label{eq:cover-cond-5}%
\end{equation}

\end{definition}

Conditions (\ref{eq:cover-cond-3}) and (\ref{eq:cover-cond-2}) mean that, in
the $y$ (stable) projection, $f_{\phi}$ contracts the set $D_{\phi}$ strictly
inside $\bar{B}_{s}^{r}$. Conditions (\ref{eq:cover-cond-4}) and
(\ref{eq:cover-cond-1}) mean that, in the $x$ (unstable) projection, $f_{\phi
}$ expands the set $D_{\phi}$ strictly outside $\bar{B}_{s}^{r}$. The final
assumption (\ref{eq:cover-cond-5}) is needed to ensure that the image of
$D_{\phi}$ by $f_{\phi}$ intersects $D_{\phi}$. Without assumption
\eqref{eq:cover-cond-5}, all other assumptions (\ref{eq:cover-cond-1}%
),\ldots,(\ref{eq:cover-cond-4}) could easily follow from having image of $D$
disjoint with $D$.

Covering relations are tools which can be used to ensure existence of an
invariant set in $D$. To prove that this set is a manifold we shall need
additional assumptions. These shall be expressed by ``cone conditions''. To
introduce these conditions, first we need some notations.

Let $Q_{h},Q_{v}:\mathbb{R}^{c}\times\mathbb{R}^{s}\times\mathbb{R}%
^{u}\rightarrow\mathbb{R}$ be functions defined by
\begin{align}
Q_{h}(\theta,x,y)  &  =-\gamma||\theta||^{2}+\alpha_{h}||x||^{2}-\beta
_{h}||y||^{2},\label{eq:Qh-def}\\
Q_{v}(\theta,x,y)  &  =-\gamma||\theta||^{2}-\alpha_{v}||x||^{2}+\beta
_{v}||y||^{2},\nonumber
\end{align}
with $\gamma,\alpha_{h},\alpha_{v},\beta_{h},\beta_{v}>0$ and%
\begin{equation}
\alpha_{h}>\alpha_{v}\quad\text{and}\quad\beta_{v}>\beta_{h}.
\label{eq:hor-vert-coeff-alignment}%
\end{equation}

\begin{figure}[ptb]
\begin{center}
\includegraphics[width=4.8in]{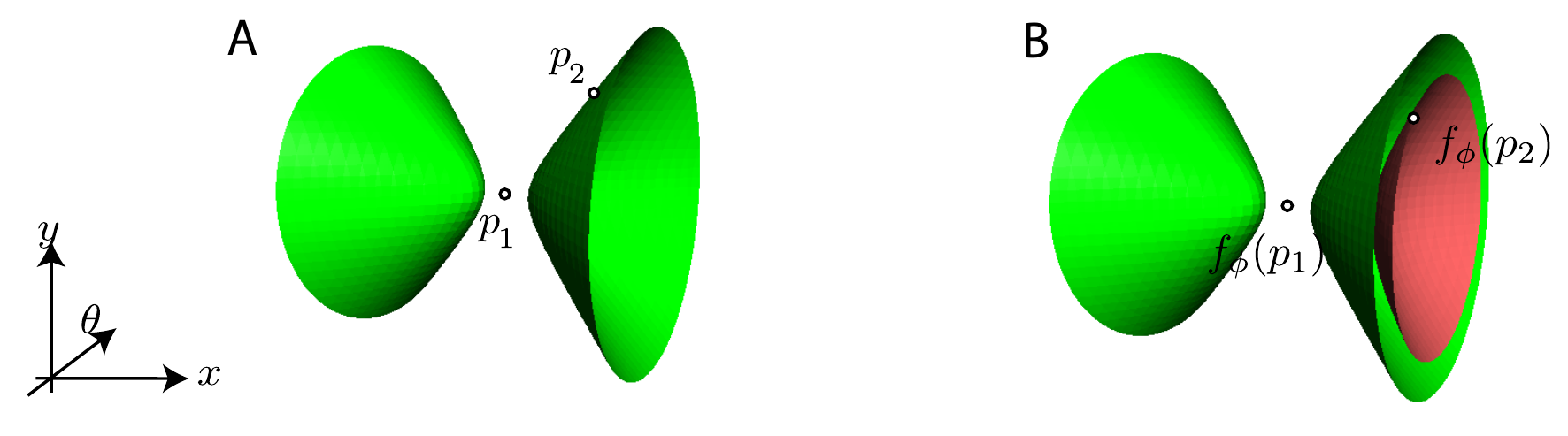}
\end{center}
\caption{An example of a function $f$, which satisfies cone conditions: A. Two
points $p_{1},p_{2}$ for which $Q_{h}(p_{1}-p_{2})=c>0$. B. Difference of the
images of the points lie on a cone $Q_{h}(f_{\phi}(p_{1})-f_{\phi}(p_{2}%
))>mc$. Similar condition (but with reversed roles of the $x$ and $y$
coordinates) needs to hold for the inverse map.}%
\label{fig:cc}%
\end{figure}

\begin{definition}
\label{def:cc} We say that a map $f:U\rightarrow U$ satisfies \emph{cone
conditions} in $D$ if there exists an $m>1$ such that

\begin{enumerate}
\item for any two points $p_{1},p_{2}\in D_{\phi}$ satisfying $p_{1}\neq
p_{2}$ and $Q_{h}(p_{1}-p_{2})\geq0$ we have%
\begin{equation}
Q_{h}(f_{\phi}(p_{1})-f_{\phi}(p_{2}))>mQ_{h}(p_{1}-p_{2}),\label{eq:cc-forw}%
\end{equation}

\item for any two points $p_{1},p_{2}\in D_{\phi}$ satisfying $p_{1}\neq
p_{2}$ and $Q_{v}(p_{1}-p_{2})\geq0$ we have%
\begin{equation}
Q_{v}(f_{\phi}^{-1}(p_{1})-f_{\phi}^{-1}(p_{2}))>mQ_{v}(p_{1}-p_{2}%
).\label{eq:cc-back}%
\end{equation}

\end{enumerate}
\end{definition}

Definition \ref{def:cc} intuitively states that if we have two points that lie
horizontally with respect to each other, then their images are going to be
pulled apart in the horizontal, $x$ coordinate (see Figure \ref{fig:cc}). If
on the other have we have two points that lie vertically with respect to each
other, then their pre-images are going to be pulled apart in the vertical, $y$
coordinate. 

We now give definitions of horizontal discs and vertical discs. These will be
 building blocks in our construction of invariant manifolds.

\begin{figure}[ptb]
\begin{center}
\includegraphics[height=1in]{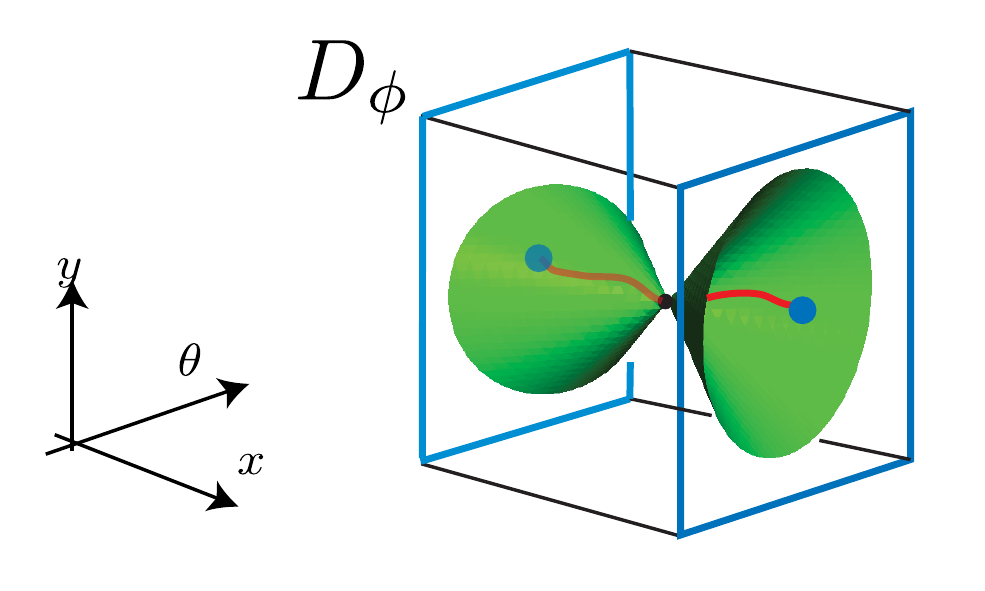}
\end{center}
\caption{A horizontal disc in $D_{\phi}$.}%
\label{fig:hor-disc}%
\end{figure}

\begin{definition}
We say that a continuous monomorphism $\mathbf{h}:\bar{B}_{u}^{r}\rightarrow
D_{\phi}$ is a \emph{horizontal disc} if $\pi_{x}\mathbf{h}(x)=x$ and for any
$x_{1},x_{2}\in\bar{B}_{u}^{r}$%
\begin{equation}
Q_{h}(\mathbf{h}(x_{1})-\mathbf{h}(x_{2}))\geq0. \label{eq:cc-hor-disc}%
\end{equation}
Thus, to any point $x$ in the graph $\mathbf{h}(x)$ we can attach a horizontal
cone, so that the graph always remains entirely inside the cone (see Figure
\ref{fig:hor-disc}).
\end{definition}

\begin{definition}
We say that a continuous monomorphism $\mathbf{v}:\bar{B}_{s}^{r}\rightarrow
D_{\phi}$ is a \emph{vertical disc} if $\pi_{y}\mathbf{v}(y)=y$ and for any
$y_{1},y_{2}\in\bar{B}_{s}^{r}$%
\[
Q_{v}(\mathbf{v}(y_{1})-\mathbf{v}(y_{2}))\geq0.
\]
Thus, to any point $y$ in the graph $\mathbf{v}(y)$ we can attach a vertical
cone, so that the graph always remains entirely inside the cone.
\end{definition}

The following lemma is a key auxiliary result for the proof of Theorem
\ref{th:nhim-th-flows}, which is the main result of this section. Roughly
speaking, it states that under appropriate conditions, an image of a
horizontal disc is a horizontal disc.

\begin{lemma}
\label{lem:image-of-hor-disc}Let $\mathbf{h}_{1}$ be a horizontal disc. If $f$
satisfies covering and cone conditions in $D$, then there exists a horizontal
disc $\mathbf{h}_{2}$ such that%
\[
\{p:\pi_{x,y}p\in \bar B_u^r \times \bar B_s^r\} \cap f_{\phi}(\mathbf{h}_{1}(\bar{B}_{u}^{r}))=\mathbf{h}_{2}(\bar{B}_{u}^{r}).
\]
Moreover, if $H(\phi^{-1}(\mathbf{h}_{1}(\bar{B}_{u}^{r})))<h$, and for any
$p\in D_{\phi}$ such that $H(\phi^{-1}(p))<h$ we have%
\begin{equation}
\pi_{\theta}(f_{\phi}(p))\in\bar{B}_{c}^{R},\label{eq:h-bound-for-maps}%
\end{equation}
then
\[
\mathbf{h}_{2}(\bar{B}_{u}^{r})\subset D_{\phi}\quad\text{and\quad}H(\phi
^{-1}(\mathbf{h}_{2}(\bar{B}_{u}^{r})))<h.
\]

\end{lemma}

\begin{proof}
Without loss of generality we can assume that $\phi$ is equal to identity.
Thus we can set $D_{\phi}=D$ and $f_{\phi}=f$.

Let $\mathbf{h}$ be any horizontal disc, then by (\ref{eq:Qh-def}),
(\ref{eq:cc-hor-disc}) and (\ref{eq:cc-forw}) for $x_{1}\neq x_{2}$%
\begin{align}
\alpha_{h}\left\Vert \pi_{x}f(\mathbf{h}(x_{1}))-\pi_{x}f(\mathbf{h}%
(x_{2}))\right\Vert ^{2} &  \geq Q_{h}(f(\mathbf{h}(x_{1}))-f(\mathbf{h}%
(x_{2})))\label{eq:cone-cond-proof1}\\
&  >mQ_{h}(\mathbf{h}(x_{1})-\mathbf{h}(x_{2}))\nonumber\\
&  \geq0,\nonumber
\end{align}
which means that $\pi_{x}\circ f\circ\mathbf{h}$ is injective.

Let us define a function $F:\bar{B}_{u}^{r}\rightarrow\mathbb{R}^{u}$ as
follows%
\[
F(x):=\pi_{x}(f(\mathbf{h}_{1}(x))).
\]
We shall first show that there exists an $x_{0}\in\bar{B}_{u}$ such that
$F(x_{0})\in\bar{B}_{u}^{r}$. Using notations $\mathbf{h}_{1}(x)=(h_{\theta
}(x),x,h_{y}(x))$ we can define a family of horizontal discs $\mathbf{h}%
_{\alpha}(x)=(\alpha h_{\theta}(x),x,\alpha h_{y}(x)).$ We define a function
$l:[0,1]\times\bar{B}_{u}^{r}\rightarrow\mathbb{R}^{u}$ as
\[
l(\alpha,x):=\pi_{x}\circ f\circ\mathbf{h}_{\alpha}(x).
\]
By (\ref{eq:Dphi-minus}) and (\ref{eq:cover-cond-1}), since $\mathbf{h}%
_{\alpha}(\partial\bar{B}_{u}^{r})\subset D^{-},$ for any $\alpha\in
\lbrack0,1]$ we have $l(\alpha,\partial\bar{B}_{u}^{r})\cap\bar{B}_{u}%
^{r}=\emptyset.$ Since, as shown at the beginning of the proof,
\[
l(\alpha,\cdot):=\pi_{x}\circ f\circ\mathbf{h}_{\alpha}:\bar{B}_{u}%
^{r}\rightarrow\mathbb{R}^{u}%
\]
is a continuous monomorphism, we either have $l(\alpha,\bar{B}_{u}^{r}%
)\cap\bar{B}_{u}^{r}=\emptyset$ or $\bar{B}_{u}^{r}\subset\mathrm{int}%
(l(\alpha,\bar{B}_{u}^{r})).$ This also means that%
\[
\inf\{\left\Vert l(\alpha,0)-x\right\Vert :x\in\partial\bar{B}_{u}^{r}\}>0,
\]
and thus the function $\delta:[0,1]\rightarrow\mathbb{R}$ defined as%
\[
\delta(\alpha):=\left\{
\begin{array}
[c]{lll}%
0 &  & l(\alpha,0)\in\bar{B}_{u}^{r}\\
1 &  & l(\alpha,0)\notin\bar{B}_{u}^{r},
\end{array}
\right.
\]
is continuous. We have%
\begin{align*}
\pi_{x}\mathbf{h}_{\alpha=0}(0) &  =0,\\
\pi_{y}\mathbf{h}_{\alpha=0}(0) &  =0,
\end{align*}
so condition (\ref{eq:cover-cond-5}) implies $l(0,0)=\pi_{x}\circ
f(\mathbf{h}_{\alpha=0}(0))\in B_{u}^{r},$ hence $\delta(0)=0.$ Suppose, to
obtain a contradiction, that $F(x)\notin\bar{B}_{u}^{r}$ for all $x\in\bar
{B}_{u}^{r}$. This would mean that in particular $F(0)=l(1,0)\notin\bar{B}%
_{u}^{r}$, hence $\delta(1)=1.$ This contradicts the fact that $\delta(0)=0$
and $\delta$ is continuous.

We have shown that there exists an $x_{0}\in\bar{B}_{u}^{r}$ such that
$F(x_{0})\in\bar{B}_{u}^{r}.$ From (\ref{eq:cover-cond-1}) follows that
$F(\partial\bar{B}_{u}^{r})\cap\bar{B}_{u}^{r}=\emptyset.$ We also know that
$F=\pi_{x}\circ f\circ\mathbf{h}_{1}$ is continuous and injective. Putting
these facts together gives $\bar{B}_{u}^{r}\subset F(\bar{B}_{u}^{r}).$ This
means that for any $v\in\bar{B}_{u}^{r}$ there exists a unique $x=x(v)\in
B_{u}^{r}$ such that $F(x)=v.$ We define
\[
\mathbf{h}_{2}(v)=(h_{2,\theta}(v),v,h_{2,y}(v)):=(\pi_{\theta}\circ
f\circ\mathbf{h}_{1}(x(v)),v,\pi_{y}\circ f\circ\mathbf{h}_{1}(x(v))).
\]
For any $v_{1}\neq v_{2}$, $v_{1},v_{2}\in\overline{B}_{u}$, by
(\ref{eq:cc-hor-disc}) and (\ref{eq:cc-forw}) we have%
\begin{align}
Q_{h}\left(  \mathbf{h}_{2}(v_{1})-\mathbf{h}_{2}(v_{2})\right)   &
=Q_{h}(f\circ\mathbf{h}_{1}(x(v_{1}))-f\circ\mathbf{h}_{1}(x(v_{2}%
)))\label{eq:cc-image-ind}\\
&  >mQ_{h}(\mathbf{h}_{1}(x(v_{1}))-\mathbf{h}_{1}(x(v_{2})))\nonumber\\
&  \geq0.\nonumber
\end{align}
Since $Q_{h}\left(  \mathbf{h}_{2}(v_{1})-\mathbf{h}_{2}(v_{2})\right)  >0,$%
\begin{align*}
\alpha_{h}\left\Vert v_{1}-v_{2}\right\Vert  &  >\beta_{h}\left\Vert
h_{2,y}(v_{1})-h_{2,y}(v_{2})\right\Vert ^{2}+\gamma\left\Vert h_{2,\theta
}(v_{1})-h_{2,\theta}(v_{2})\right\Vert ^{2}\\
&  \geq\min(\beta_{h},\gamma)\left\Vert (h_{2,\theta},h_{2,y})(v_{1}%
)-(h_{2,\theta},h_{2,y})(v_{2})\right\Vert ^{2},
\end{align*}
and therefore $\mathbf{h}_{2}$ is continuous.

Finally let us note that (\ref{eq:map-integral}) and $\mathbf{h}_{2}%
(v)=f\circ\mathbf{h}_{1}(x(v))$ implies $H(\mathbf{h}_{2}(\bar{B}_{u}%
^{r}))=H(\mathbf{h}_{1}(\bar{B}_{u}^{r}))<h.$ This by
(\ref{eq:h-bound-for-maps}) implies that $\mathbf{h}_{2}(\bar{B}_{u}%
^{r})\subset D$.
\end{proof}

Next lemma follows from mirror arguments.

\begin{lemma}
\label{lem:image-of-ver-disc}Let $\mathbf{v}_{1}$ be a vertical disc. If $f$
satisfies covering and cone conditions in $D$, then there exists a vertical
disc $\mathbf{v}_{2}$ such that%
\[
\{p:\pi_{x,y}p\in \bar B_u^r \times \bar B_s^r\} \cap f_{\phi}(\mathbf{v}_{1}(\bar{B}_{s}^{r}))=\mathbf{v}_{2}(\bar{B}_{s}^{r}).
\]
Moreover, if $H(\phi^{-1}(\mathbf{v}_{1}(\bar{B}_{s}^{r})))<h$, and for any
$p\in D_{\phi}$ such that $H(\phi^{-1}(p))<h$ we have%
\begin{equation}
\pi_{\theta}(f_{\phi}^{-1}(p))\in\bar{B}_{c}^{R},\label{eq:v-bound-for-maps}%
\end{equation}
then
\[
\mathbf{v}_{2}(\bar{B}_{s}^{r})\subset D_{\phi}\quad\text{and\quad}H(\phi
^{-1}(\mathbf{v}_{2}(\bar{B}_{s}^{r})))<h.
\]

\end{lemma}

We are now ready to state our main result for maps, which will be the main
tool for the proof of Theorem \ref{th:main}.

\begin{theorem}
\label{th:nhim-th-flows} If $f$ satisfies covering and cone conditions in $D$,
and in addition for any $p\in D_{\phi}$ with $H(\phi^{-1}(p))<h$ we have%
\begin{equation}
\pi_{\theta}f_{\phi}(p)\in\bar{B}_{c}^{R}\quad\text{and}\quad\pi_{\theta
}f_{\phi}^{-1}(p)\in\bar{B}_{c}^{R},\label{eq:map-cond-on-centre}%
\end{equation}
then there exists a $C^{0}$ function $\chi:\bar{B}_{c}^{R-v}\rightarrow
D_{\phi}$ such that

\begin{enumerate}
\item For any $\theta\in\bar{B}_{c}^{R-v}$ we have $\pi_{\theta}\chi
(\theta)=\theta$ and
\[
f^{n}(\phi^{-1}(\chi(\theta)))\in D\quad\text{for all }n\in\mathbb{Z}.
\]

\item If for some $p\in\phi^{-1}(\bar{B}_{c}^{R-v}\times\bar{B}_{u}^{r}%
\times\bar{B}_{s}^{r})$ we have
\[
f^{n}(p)\in D\quad\text{for all }n\in\mathbb{Z,}%
\]
then there exists a $\theta\in\bar{B}_{c}^{R-v}$ such that $p=\phi^{-1}%
(\chi(\theta))$.
\end{enumerate}
\end{theorem}

\begin{proof}
Without loss of generality we assume that $\phi$ is equal to identity,
which means that $D_{\phi}=D$ and $f_{\phi}=f$.

Let $\theta_{0}\in\bar{B}_{c}^{R-v}$ and $y_{0}\in\bar{B}_{s}^{r}$. Let
$\mathbf{h}_{1}:\bar{B}_{u}^{r}\rightarrow D$ be a horizontal disc defined by%
\[
\mathbf{h}_{1}(x)=(\theta_{0},x,y_{0}).
\]
Clearly $\mathbf{h}_{1}$ satisfies cone conditions and also by
(\ref{eq:h-boundary-bound}), $H(\phi^{-1}(\mathbf{h}_{1}(\bar{B}_{u}^{r}%
)))<h$. Applying inductively Lemma \ref{lem:image-of-hor-disc} we obtain a
sequence of horizontal discs $\mathbf{h}_{1},\mathbf{h}_{2},\ldots$ such that%
\[
f(\mathbf{h}_{i-1}(\bar{B}_{u}^{r}))\cap D=\mathbf{h}_{i}(\bar{B}_{u}%
^{r})\quad\text{and\quad}H(\mathbf{h}_{i}(\bar{B}_{u}^{r}))<h.
\]
This by compactness of $\bar{B}_{u}^{r}$ ensures existence of a point
$x_{0}^{\ast}\in\bar{B}_{u}^{r}$ such that for all $n\in\mathbb{N}$%
\begin{equation}
f^{n}(\mathbf{h}_{1}(x_{0}^{\ast}))\in D.\label{eq:map-temp1}%
\end{equation}
Suppose that we have two points $x_{0}^{1}$ and $x_{0}^{2}$ which satisfy
(\ref{eq:map-temp1}). Then by (\ref{eq:cc-image-ind}) we have%
\begin{align}
\alpha_{h}r^{2} &  \geq\alpha_{h}\left\Vert \pi_{x}\left(  f^{n}%
(\mathbf{h}_{1}(x_{0}^{1}))-f^{n}(\mathbf{h}_{1}(x_{0}^{2}))\right)
\right\Vert ^{2}\nonumber\\
&  >Q_{h}\left(  f^{n}(\mathbf{h}_{1}(x_{0}^{1}))-f^{n}(\mathbf{h}_{1}%
(x_{0}^{2}))\right)  \nonumber\\
&  >mQ_{h}\left(  f^{n-1}(\mathbf{h}_{1}(x_{0}^{1}))-f^{n-1}(\mathbf{h}%
_{1}(x_{0}^{2}))\right)  \label{eq:hor-induction-expansion}\\
&  \ldots\nonumber\\
&  >m^{n}Q_{h}\left(  \mathbf{h}_{1}(x_{0}^{1})-\mathbf{h}_{1}(x_{0}%
^{2})\right)  ,\nonumber
\end{align}
which since $m>1$ cannot hold for all $n$. This means that functions
$W^{cs}:\bar{B}_{c}^{R-v}\times\bar{B}_{s}^{r}\rightarrow D$, $w^{cs}:\bar
{B}_{c}^{R-v}\times\bar{B}_{s}^{r}\rightarrow\bar{B}_{u}^{r}$ given as
\[
W^{cs}(\theta_{0},y_{0})=(\theta_{0},w^{cs}(\theta_{0},y_{0}),y_{0}%
):=(\theta_{0},x_{0}^{\ast},y_{0}),
\]
are properly defined. Note that by a similar argument to
(\ref{eq:hor-induction-expansion}), for any $(\theta_{1},y_{1})\neq(\theta
_{2},y_{2})$ we must have%
\begin{equation}
0>Q_{h}(W^{cs}(\theta_{2},y_{2})-W^{cs}(\theta_{1},y_{1}%
)).\label{eq:Ws-in-negative-cone}%
\end{equation}
This gives%
\begin{align*}
\max(\gamma,\beta_{h})\left\Vert (\theta_{2},y_{2})-(\theta_{1},y_{1}%
)\right\Vert ^{2} &  \geq\gamma||\theta_{2}-\theta_{1}||^{2}+\beta_{h}%
||y_{2}-y_{1}||^{2}\\
&  >\alpha_{h}||w^{cs}(\theta_{2},y_{2})-w^{cs}(\theta_{2},y_{2})||^{2},
\end{align*}
which means that $w^{cs}$ is Lipschitz with a constant
\begin{equation}
L_{s}=\sqrt{\frac{\max(\gamma,\beta_{h})}{\alpha_{h}}}.\label{eq:lip-ws}%
\end{equation}

Mirror arguments, involving Lemma \ref{lem:image-of-ver-disc}, give existence of
functions $W^{cu}:\bar{B}_{c}^{R-v}\times\bar{B}_{u}^{r}\rightarrow D$,
$w^{cu}:\bar{B}_{c}^{R-v}\times\bar{B}_{u}^{r}\rightarrow\bar{B}_{s}^{r}$
\[
W^{cu}(\theta,x)=(\theta,x,w^{cu}(\theta,y)),
\]
such that for any point $(\theta,x)\in\bar{B}_{c}^{R-v}\times\bar{B}_{u}^{r}$
and all $n\in\mathbb{N}$%
\[
f^{-n}(W^{cu}(\theta,x))\in D.
\]
Also $w^{cu}$ is Lipschitz with a constant
\begin{equation}
L_{u}=\sqrt{\frac{\max(\gamma,\alpha_{v})}{\beta_{v}}}.\label{eq:lip-wu}%
\end{equation}

We shall show that for any $\theta\in\bar{B}_{c}^{R-v}$ the sets
$W^{cs}(\theta,\bar{B}_{s}^{r})$ and $W^{cu}(\theta,\bar{B}_{u}^{r})$
intersect. Let us define $P_{\theta}:\bar{B}_{u}^{r}\times\bar{B}_{s}%
^{r}\rightarrow\bar{B}_{u}^{r}\times\bar{B}_{s}^{r}$ as%
\[
P_{\theta}(x,y):=\left(  \pi_{x}W^{cs}(\theta,y),\pi_{y}W^{cu}(\theta
,x)\right)  .
\]
Since $P_{\theta}$ is continuous, by the Brouwer fixed point theorem 
 there exists an $(x_{0},y_{0})$ such that $P_{\theta}(x_{0}%
,y_{0})=\left(  x_{0},y_{0}\right)  .$ This means that%
\[
W^{cs}(\theta,y_{0})=\left(  \theta,w^{cs}(\theta,y_{0}),y_{0}\right)
=\left(  \theta,x_{0},w^{cu}(\theta,y_{0})\right)  =W^{cu}(\theta,x_{0}).
\]

Now we shall show that for any given $\theta\in\bar{B}_{c}^{R-v}$ there exists
only a single point of such intersection. Suppose that for some $\theta\in
\bar{B}_{c}^{R-v}$ there exist $\left(  x_{1},y_{1}\right)  ,\left(
x_{2},y_{2}\right)  \in\bar{B}_{u}^{r}\times\bar{B}_{s}^{r}$, $\left(
x_{1},y_{1}\right)  \neq\left(  x_{2},y_{2}\right)  $ such that
\[
W^{cs}(\theta,y_{1})=W^{cu}(\theta,x_{1})\quad\text{and\quad}W^{cs}%
(\theta,y_{2})=W^{cu}(\theta,x_{2}).
\]
We would then have $W^{cs}(\theta,y_{m})=W^{cu}(\theta,x_{m})=(\theta
,x_{m},y_{m})$ for $m=1,2.$

From (\ref{eq:Ws-in-negative-cone}) follows that%
\[
0>Q_{h}\left(  W^{cs}(\theta,y_{1})-W^{cs}(\theta,y_{2})\right)  =Q_{h}\left(
(\theta,x_{1},y_{1})-(\theta,x_{2},y_{2})\right)  ,
\]
and by mirror argument%
\[
0>Q_{v}\left(  W^{cu}(x_{1},\lambda)-W^{cu}(x_{2},\lambda)\right)
=Q_{v}\left(  (\theta,x_{1},y_{1})-(\theta,x_{2},y_{2})\right)
\]
which implies that
\[
0>(\alpha_{h}-\alpha_{v})\left\Vert x_{1}-x_{2}\right\Vert ^{2}+\left(
\beta_{v}-\beta_{h}\right)  \left\Vert y_{1}-y_{2}\right\Vert ^{2},
\]
which contradicts (\ref{eq:hor-vert-coeff-alignment}).

We now define $\chi(\theta)=(\theta,\chi_{x,y}(\theta)):=(\theta,x_{0},y_{0})$
for $x_{0}=x_{0}(\theta),$ $y_{0}=y_{0}(\theta)$ such that $W^{cs}%
(\theta,y_{0})=W^{cu}(\theta,x_{0}).$ By previous arguments we know that
$\chi$ is properly defined. We need to show continuity. Let us take any $\theta_{1},\theta_{2}\in\bar{B}_{c}^{R-v}$. From
(\ref{eq:Ws-in-negative-cone}) follows that
\begin{equation}
Q_{h}\left(  \chi(\theta_{1})-\chi(\theta_{2})\right)  =Q_{h}\left(
W^{cs}(\theta_{1},y_{0}(\theta_{1}))-W^{cs}(\theta_{2},y_{0}\left(  \theta
_{2}\right)  )\right)  <0,\label{eq:cc-for-chi-h}%
\end{equation}
and by mirror argument%
\begin{equation}
Q_{v}\left(  \chi(\theta_{1})-\chi(\theta_{2})\right)  =Q_{v}\left(
W^{cu}(\theta_{1},x_{0}(\theta_{1}))-W^{cu}(\theta_{2},x_{0}\left(  \theta
_{2}\right)  )\right)  <0.\label{eq:cc-for-chi-v}%
\end{equation}
From (\ref{eq:cc-for-chi-h}), (\ref{eq:cc-for-chi-v}) follows that
\begin{align*}
\alpha_{h}\left\Vert x_{0}\left(  \theta_{1}\right)  -x_{0}\left(  \theta
_{2}\right)  \right\Vert ^{2}-\beta_{h}\left\Vert y_{0}\left(  \theta
_{1}\right)  -y_{0}\left(  \theta_{2}\right)  \right\Vert ^{2} &
<\gamma\left\Vert \theta_{1}-\theta_{2}\right\Vert ^{2},\\
-\alpha_{v}\left\Vert x_{0}\left(  \theta_{1}\right)  -x_{0}\left(  \theta
_{2}\right)  \right\Vert ^{2}+\beta_{v}\left\Vert y_{0}\left(  \theta
_{1}\right)  -y_{0}\left(  \theta_{2}\right)  \right\Vert ^{2} &
<\gamma\left\Vert \theta_{1}-\theta_{2}\right\Vert ^{2},
\end{align*}%
\[
\left(  \alpha_{h}-\alpha_{v}\right)  \left\Vert x_{0}\left(  \theta
_{1}\right)  -x_{0}\left(  \theta_{2}\right)  \right\Vert ^{2}+\left(
\beta_{v}-\beta_{h}\right)  \left\Vert y_{0}\left(  \theta_{1}\right)
-y_{0}\left(  \theta_{2}\right)  \right\Vert ^{2}<2\gamma\left\Vert \theta
_{1}-\theta_{2}\right\Vert ^{2},
\]
which gives
\begin{equation}
\left\Vert \chi_{x,y}\left(  \theta_{1}\right)  -\chi_{x,y}\left(  \theta
_{2}\right)  \right\Vert <\sqrt{\frac{2\gamma}{\min(\alpha_{h}-\alpha
_{v},\beta_{v}-\beta_{h})}}\left\Vert \theta_{1}-\theta_{2}\right\Vert
,\label{eq:chi-Lip-bound}%
\end{equation}
and by (\ref{eq:hor-vert-coeff-alignment}) implies Lipschitz bounds for
$\chi_{x,y}$ and continuity of $\chi$.
\end{proof}



\section{Proof of the main theorem\label{sec:proof-of-main}}

In this section we shall show that assumptions of Theorem \ref{th:main} imply
that a map induced as a shift along a trajectory of the flow of
(\ref{eq:our-ode}) for sufficiently small time satisfies covering and cone
conditions. This will allow us to apply Theorem \ref{th:nhim-th-flows} to
prove Theorem \ref{th:main}.

We start with a lemma which shows that assumptions of Theorem \ref{th:main}
imply covering conditions for a shift along the trajectory of (\ref{eq:our-ode}).

\begin{lemma}
\label{lem:covering-for-flows}Assume that for any $p\in\bar{B}_{c}^{R-v}%
\times\{0\}\times\{0\}$ assumptions (\ref{eq:mth-as-cover-1}),
(\ref{eq:mth-as-cover-2}), (\ref{eq:mth-as-cover-3}) of Theorem \ref{th:main}
hold, then for sufficiently small $\tau>0$ and all $t\in(0,\tau]$ a function%
\[
f(\mathbf{x}):=\Phi(t,\mathbf{x})
\]
satisfies covering conditions.
\end{lemma}

\begin{proof}
Without loss of generality we assume that $\phi=id$. Let $q\in B_{c}^{R}%
\times\{0\}\times\{0\}\cap N_{p}$. By (\ref{eq:mth-as-cover-3}), for
sufficiently small $t$%
\begin{align}
\left\Vert \pi_{x}\Phi(t,q)\right\Vert  &  =\left\Vert \pi_{x}\left(
\Phi(0,q)+\frac{d}{dt}\Phi(0,q)t+o(t)\right)  \right\Vert \nonumber\\
&  =\left\Vert 0+t\pi_{x}F(q)+o(t)\right\Vert \nonumber\\
&  <\left\vert t\right\vert rE_{u}. \label{eq:Fx-on-central}%
\end{align}
Analogous computation yields%
\begin{equation}
\left\Vert \pi_{y}\Phi(t,q)\right\Vert <\left\vert t\right\vert rE_{s}.
\label{eq:Fy-on-central}%
\end{equation}

In later parts of the proof we shall use the fact that for any $q_{1},q_{2}\in
N_{p}$%
\begin{equation}
F(q_{1})-F(q_{2})=\int_{0}^{1}dF\left(  q_{2}+s\left(  q_{1}-q_{2}\right)
\right)  ds\left(  q_{1}-q_{2}\right)  . \label{eq:difference-from-derivative}%
\end{equation}

Now we shall prove (\ref{eq:cover-cond-1}). Let $q=(\theta,x,y)\in D_{\phi
}^{-}\cap N_{p},$ which means that $\left\Vert x\right\Vert =r.$ Using
$\frac{d}{dt}\Phi(t,q)|_{t=0}=F(q),$ $\Phi(0,q)=q,$ and
(\ref{eq:difference-from-derivative}) we have
\begin{align*}
&  \frac{d}{dt}\left\Vert \pi_{x}\left(  \Phi(t,q)-\Phi(t,\left(
\theta,0,0\right)  )\right)  \right\Vert ^{2}|_{t=0}\\
&  =\frac{d}{dt} \left. \left( \pi_{x}\left(  \Phi(t,q)-\Phi(t,\left(  \theta,0,0\right)
)\right)  ^{T}\pi_{x}\left(  \Phi(t,q)-\Phi(t,\left(  \theta,0,0\right)
)\right) \right) \right|_{t=0}\\
&  =2\pi_{x}\left(  q-\left(  \theta,0,0\right)  \right)  ^{T}\pi_{x}\left(
F(q)-F\left(  \theta,0,0\right)  \right) \\
&  =2x^{T}\pi_{x}\left(  \int_{0}^{1}dF\left(  \theta,sx,sy\right)  ds\left(
0,x,y\right)  \right) \\
&  =2x^{T}\left(  Ax+e_{u}y\right)  ,
\end{align*}
where%
\[
A=\int_{0}^{1}\frac{\partial(\pi_{x}F)}{\partial x}\left(  \theta
,sx,sy\right)  ds,\quad\quad e_{u}=\int_{0}^{1}\frac{\partial\left(  \pi
_{x}F\right)  }{\partial y}\left(  \theta,sx,sy\right)  ds.
\]
From bounds (\ref{eq:est-expansion}) and (\ref{eq:est-epsilons}) we thus
obtain%
\begin{equation}
\frac{d}{dt}\left\Vert \pi_{x}\left(  \Phi(t,q)-\Phi(t,\left(  \theta
,0,0\right)  )\right)  \right\Vert ^{2}|_{t=0}>2\left(  r^{2}\delta
^{u}-\left\Vert x\right\Vert \left\Vert e_{u}\right\Vert \left\Vert
y\right\Vert \right)  >2r^{2}\left(  \delta^{u}-\varepsilon_{u}\right)  .
\label{eq:dt-on-Phix}%
\end{equation}
Using the same arguments we can also show that for any $q=(\theta,x,y)\in
N_{p}$%
\begin{equation}
\frac{d}{dt}\left\Vert \pi_{y}\left(  \Phi(t,q)-\Phi(t,\left(  \theta
,0,0\right)  )\right)  \right\Vert ^{2}|_{t=0}<2\left\Vert y\right\Vert
\left(  \varepsilon_{s}r-\left\Vert y\right\Vert \delta^{s}\right)
\label{eq:dt-on-Phiy}%
\end{equation}

Combining (\ref{eq:Fx-on-central}) (\ref{eq:dt-on-Phix}) and
(\ref{eq:mth-as-cover-1}), for sufficiently small $t>0$ gives%
\begin{align}
\left\Vert \pi_{x}f(q)\right\Vert  &  =\left\Vert \pi_{x}\Phi(t,q)\right\Vert
\nonumber\\
&  \geq\left\Vert \pi_{x}\left(  \Phi(t,q)-\Phi(t,(\theta,0,0))\right)
\right\Vert -\left\Vert \pi_{x}\Phi(t,(\theta,0,0))\right\Vert \nonumber\\
&  >\sqrt{\left\Vert \pi_{x}\left(  \Phi(t,q)-\Phi(t,(\theta,0,0))\right)
\right\Vert ^{2}}-trE_{u}\label{eq:fx-bound-r}\\
&  >\sqrt{\left\Vert \pi_{x}\left(  \Phi(0,q)-\Phi(0,(\theta,0,0))\right)
\right\Vert ^{2}+t2\left(  \delta^{u}-\varepsilon_{u}r^{2}\right)  }%
-trE_{u}\nonumber\\
&  =\sqrt{r^{2}+t2r^{2}\left(  \delta^{u}-\varepsilon_{u}\right)  }%
-trE_{u}\nonumber\\
&  >r.\nonumber
\end{align}
This establishes (\ref{eq:cover-cond-1}). Now we shall show
(\ref{eq:cover-cond-3}). For any $q=(\theta,x,y)\in D_{\phi}$ and sufficiently
small $t>0,$ analogous derivation to (\ref{eq:fx-bound-r}) (for these
computations we use estimates (\ref{eq:Fy-on-central}), (\ref{eq:dt-on-Phiy}))
give%
\begin{equation}
\left\Vert \pi_{y}f(q)\right\Vert <\sqrt{\left\Vert y\right\Vert
^{2}+t2\left\Vert y\right\Vert \left(  \varepsilon_{s}r-\left\Vert
y\right\Vert \delta^{s}\right)  }+trE_{s}. \label{eq:fy-bound-r}%
\end{equation}
Since $\left\Vert y\right\Vert \leq r$ by (\ref{eq:mth-as-cover-2}), for
sufficiently small $t>0,$ inequality (\ref{eq:fy-bound-r}) implies that
$\left\Vert \pi_{y}f(q)\right\Vert <r$ and hence establishes
(\ref{eq:cover-cond-3}).

Proof of (\ref{eq:cover-cond-2}) and (\ref{eq:cover-cond-4}) follows from
analogous arguments with $t<0$.

Conditions (\ref{eq:cover-cond-5}) hold for sufficiently small $t.$ This
follows from continuity of $\Phi(p,t)$ with respect to $t$ since
\[
f(p)=\Phi(t,p)\qquad f^{-1}(p)=\Phi(-t,p),
\]
and for $p\in\bar{B}_{c}^{R}\times\{0\}\times\{0\}$%
\[
\pi_{(x,y)}\Phi(0,p)=\pi_{(x,y)}p=(0,0)\in\mathrm{int}\left(  \bar{B}_{u}%
^{r}\times\bar{B}_{s}^{r}\right)  .
\]

\end{proof}

Now we shall show that assumptions of Theorem \ref{th:main} imply cone
conditions for a shift along trajectory of (\ref{eq:our-ode}). Let us start
with a simple technical lemma.

\begin{lemma}
\label{lem:C-est}Let $C=\left(  C_{ij}\right)  _{i,j=1,2,3}$ be a $\left(
c+u+s\right)  \times\left(  c+u+s\right)  $ matrix. Assume that for
$a_{i},b_{i}\in\mathbb{R}$, $i=1,2,3$ we have
\begin{align}
\inf\{x_{i}^{T}C_{ii}x_{i}  &  :\,\left\Vert x_{i}\right\Vert =1\}\geq
a_{i},\quad\text{for }i=1,2,3,\label{eq:lem-ai-est}\\
\sup\{x_{i}^{T}C_{ii}x_{i}  &  :\,\left\Vert x_{i}\right\Vert =1\}\leq
b_{i}\quad\text{for }i=1,2,3, \label{eq:lem-bi-est}%
\end{align}
then for any $x=\left(  x_{1},x_{2},x_{3}\right)  \in\mathbb{R}^{c+u+s}$%
\begin{align}
x^{T}Cx  &  \geq\left(  a_{1}-c_{1}\right)  \left\Vert x_{1}\right\Vert
^{2}+\left(  a_{2}-c_{2}\right)  \left\Vert x_{2}\right\Vert ^{2}+\left(
a_{3}-c_{3}\right)  \left\Vert x_{3}\right\Vert ^{2},\label{eq:lem-res1n}\\
x^{T}Cx  &  \leq\left(  b_{1}+c_{1}\right)  \left\Vert x_{1}\right\Vert
^{2}+\left(  b_{2}+c_{2}\right)  \left\Vert x_{2}\right\Vert ^{2}+\left(
b_{3}+c_{3}\right)  \left\Vert x_{3}\right\Vert ^{2}, \label{eq:lem-res2n}%
\end{align}
where%
\begin{align*}
c_{1}  &  =\frac{1}{2}\left(  \left\Vert C_{21}\right\Vert +\left\Vert
C_{31}\right\Vert +\left\Vert C_{12}\right\Vert +\left\Vert C_{13}\right\Vert
\right)  ,\\
c_{2}  &  =\frac{1}{2}\left(  \left\Vert C_{21}\right\Vert +\left\Vert
C_{23}\right\Vert +\left\Vert C_{12}\right\Vert +\left\Vert C_{32}\right\Vert
\right)  ,\\
c_{3}  &  =\frac{1}{2}\left(  \left\Vert C_{31}\right\Vert +\left\Vert
C_{23}\right\Vert +\left\Vert C_{13}\right\Vert +\left\Vert C_{32}\right\Vert
\right)  .
\end{align*}

\end{lemma}

\begin{proof}
The estimate (\ref{eq:lem-res1n}) follows by direct computation from
(\ref{eq:lem-ai-est}) and the fact that for any $i,j$%
\[
\pm2q_{j}^{T}C_{ji}q_{i}\geq-||C_{ji}||\left(  ||q_{j}||^{2}+||q_{i}%
||^{2}\right)  .
\]
Similarly (\ref{eq:lem-res2n}) follows from (\ref{eq:lem-bi-est}) and%
\[
\pm2q_{j}^{T}C_{ji}q_{i}\leq||C_{ji}||\left(  ||q_{j}||^{2}+||q_{i}%
||^{2}\right)  .
\]

\end{proof}

Let $I_{k}$ denote a $k\times k$ identity matrix. Let
\begin{align*}
Q_{1} &  =\mathrm{diag}(-\gamma I_{c},\alpha_{h}I_{u},-\beta_{h}I_{s}),\\
Q_{2} &  =\mathrm{diag}(-\gamma I_{c},-\alpha_{v}I_{u},\beta_{v}I_{s}),
\end{align*}
be matrices associated with quadratic forms $Q_{h}$ and $Q_{v}$ respectively.
Now we are ready to prove that assumptions of Theorem \ref{th:main} imply cone
conditions for a time shift along a trajectory map.

\begin{lemma}
\label{lem:cc-for-flows}Assume that for any $p\in\bar{B}_{c}^{R-v}%
\times\{0\}\times\{0\}$ assumption (\ref{eq:mth-as-cone1}) of Theorem
\ref{th:main} holds, then for sufficiently small $\tau>0$ and all $t\in
(0,\tau]$ a function%
\[
f(\mathbf{x}):=\Phi(t,\mathbf{x})
\]
satisfies cone conditions with a coefficient $m=1+th,$ with some constant
$h>0.$.
\end{lemma}

\begin{proof}
Let $p_{1},p_{2}\in D_{\phi}$ be such that $p_{i}=\left(  \theta_{i}%
,x_{i},y_{i}\right)  $ for $i=1,2,$ $p_{1}\neq p_{2}$ and $Q_{h}(p_{1}%
-p_{2})\geq0$. Let $p=\left(  \theta_{1},0,0\right)  \in\bar{B}_{c}^{R}%
\times\{0\}\times\{0\}$. Condition (\ref{eq:rho-cond-1}) implies that
$p_{1},p_{2}\in N_{p}$. We compute
\begin{align}
&  \frac{d}{dt}\left(  (\Phi(t,p_{1})-\Phi(t,p_{2}))^{T}Q_{1}(\Phi
(t,p_{1})-\Phi(t,p_{2}))\right)  |_{t=0}\nonumber\\
&  =2(p_{1}-p_{2})^{T}Q_{1}(F(p_{1})-F(p_{2}))\label{eq:cone-proof-dt}\\
&  =2(p_{1}-p_{2})^{T}Q_{1}B(p_{1}-p_{2}),\nonumber
\end{align}
where
\[
B=\int_{0}^{1}dF(p_{2}+t(p_{1}-p_{2}))dt\in\lbrack dF(N_{p})].
\]
For $C=Q_{1}B,$ from (\ref{eq:est-central}), (\ref{eq:est-expansion}),
(\ref{eq:est-contraction}) we have%
\begin{align}
\inf\{x_{1}^{T}C_{11}x_{1}  &  :\,\left\Vert x_{1}\right\Vert =1\}\geq-\gamma
c^{u},\nonumber\\
\inf\{x_{2}^{T}C_{22}x_{2}  &  :\,\left\Vert x_{2}\right\Vert =1\}\geq
\alpha_{h}\delta^{u},\label{eq:tempC}\\
\inf\{x_{3}^{T}C_{33}x_{3}  &  :\,\left\Vert x_{3}\right\Vert =1\}\geq
\beta_{h}\delta^{s}.\nonumber
\end{align}
Using (\ref{eq:lem-res1n}) from Lemma \ref{lem:C-est} with (\ref{eq:tempC})
and (\ref{eq:est-epsilons}), for $\kappa_{c}^{\text{forw}},\kappa
_{u}^{\text{forw}},\kappa_{s}^{\text{forw}}$ given by (\ref{eq:kappa-h}) and
$\mu_{1}\in(\max(\kappa_{c}^{\text{forw}},\kappa_{s}^{\text{forw}}),\kappa
_{u}^{\text{forw}})$ we have
\begin{align}
x^{T}Cx  &  \geq-\kappa_{c}^{\text{forw}}\gamma\left\Vert x_{1}\right\Vert
^{2}+\kappa_{u}^{\text{forw}}\alpha_{h}\left\Vert x_{2}\right\Vert ^{2}%
-\kappa_{s}^{\text{forw}}\beta_{h}\left\Vert x_{3}\right\Vert ^{2}\nonumber\\
&  >\mu_{1}\left(  -\gamma\left\Vert x_{1}\right\Vert ^{2}+\alpha
_{h}\left\Vert x_{2}\right\Vert ^{2}-\beta_{h}\left\Vert x_{3}\right\Vert
^{2}\right) \nonumber\\
&  =\mu_{1}x^{T}Q_{1}x. \label{eq:xCx-bound}%
\end{align}
The constant $\mu_{1}\in(\max(\kappa_{c}^{\text{forw}},\kappa_{s}%
^{\text{forw}}),\kappa_{u}^{\text{forw}})$ can be chosen to be greater than
zero thanks to assumption (\ref{eq:mth-as-cone1}). This means that by
(\ref{eq:cone-proof-dt}) and (\ref{eq:xCx-bound})%
\[
\frac{d}{dt}\left(  (\Phi(t,p_{1})-\Phi(t,p_{2}))^{T}Q_{1}(\Phi(t,p_{1}%
)-\Phi(t,p_{2}))\right)  |_{t=0}>2\mu_{2}Q_{h}(p_{1}-p_{2}).
\]
For sufficiently small $\tau>0$ and $t\in(0,\tau)$ we therefore have%
\begin{align*}
Q_{h}(f(p_{1})-f(p_{2}))  &  =Q_{h}(\Phi(t,p_{1})-\Phi(t,p_{2}))\\
&  =Q_{h}(p_{1}-p_{2})+t\frac{d}{dt}Q_{h}(\Phi(t,p_{1})-\Phi(t,p_{2}%
))|_{t=0}+o(t)\\
&  >(1+t2\mu_{1})Q_{h}(p_{1}-p_{2}),
\end{align*}
which establishes (\ref{eq:cc-forw}) with $m=1+t2\mu_{1}>1$.

The proof of (\ref{eq:cc-back}) is obtained analogously with $m=1+t2\mu_{2}>1$
for some $\mu_{2}<0,$ $\mu_{2}\in\left(  \kappa_{s}^{\text{back}},\min
(\kappa_{c}^{\text{back}},\kappa_{u}^{\text{back}})\right)  ,$ with negative
time $t<0$.

So far the entire argument was done for points in $N_{p}$. We can choose
$h_{p}=\min\{2\left\vert \mu_{1}\right\vert ,2\left\vert \mu_{2}\right\vert
\}$ so that (\ref{eq:cc-forw}) and (\ref{eq:cc-back}) hold for any $p_{1},p_{2}\in
N_{p}$ with a constant $m=1+\left\vert t\right\vert h_{p}$. By compactness of
$D_{\phi}$ we can now choose a $h>0$ such that (\ref{eq:cc-forw}) and
(\ref{eq:cc-back}) hold with a constant $m=1+\left\vert t\right\vert h$ for all
$p_{1},p_{2}\in D_{\phi}.$
\end{proof}

We are now ready for the proof of our main result.

\begin{proof}
[Proof of Theorem \ref{th:main}]By Lemmas \ref{lem:covering-for-flows} and
\ref{lem:cc-for-flows} we know that assumptions of Theorem \ref{th:main}
imply cone and covering conditions for a map induced by the flow by a small
time shift. Now we just need to show that for a map
\[
f(\mathbf{x}):=\Phi(t,\mathbf{x}),
\]
with sufficiently small $t>0,$ for any $p\in D_{\phi}$ with $H(\phi
^{-1}(p))<h$ we have (\ref{eq:map-cond-on-centre}). This follows from
(\ref{eq:h-boundary-bound}) and continuity of $\Phi(t,\mathbf{x})$ with
respect to $t.$ The claim now follows from Theorem \ref{th:nhim-th-flows}.
\end{proof}

By applying Theorem \ref{th:nhim-th-flows} in our proof of Theorem
\ref{th:main} we have established more than just continuity of our center
manifold. We have also obtained existence of its stable and unstable
manifolds, together with explicit Lipschitz type bounds on their slopes. This
is summarised in the following corollary.

\begin{corollary}
\label{cor:lip-bounds}During the course of the proof of Theorem
\ref{th:nhim-th-flows} we have shown that in local coordinates given by
$\phi$ the stable, unstable and center manifolds obtained by our argument are
given in terms of functions%
\begin{align*}
W^{cs}  &  :\bar{B}_{c}^{R-v}\times\bar{B}_{s}^{r}\rightarrow D_{\phi},\\
W^{cu}  &  :\bar{B}_{c}^{R-v}\times\bar{B}_{u}^{r}\rightarrow D_{\phi},\\
\chi &  :\bar{B}_{c}^{R-v}\rightarrow D_{\phi},
\end{align*}
respectively. We have also shown that these functions are of the form%
\begin{align*}
W^{cs}(\theta,y)  &  =(\theta,w^{cs}(\theta,y),y),\\
W^{cu}(\theta,x)  &  =(\theta,x,w^{cu}(\theta_{0},y)),\\
\chi(\theta)  &  =(\theta,\chi_{x,y}(\theta)),
\end{align*}
with functions $w^{cs}:\bar{B}_{c}^{R-v}\times\bar{B}_{s}^{r}\rightarrow
\bar{B}_{u}^{r}$, $w^{cu}:\bar{B}_{c}^{R-v}\times\bar{B}_{u}^{r}%
\rightarrow\bar{B}_{s}^{r}$ and $\chi_{x,y}:\bar{B}_{c}^{R-v}\rightarrow
\bar{B}_{u}^{r}\times\bar{B}_{s}^{r}$ by (\ref{eq:lip-ws}), (\ref{eq:lip-wu})
and (\ref{eq:chi-Lip-bound}) satisfying Lipschitz conditions with constants
\begin{align*}
L_{s}  &  =\sqrt{\frac{\max(\gamma,\beta_{h})}{\alpha_{h}}},\\
L_{u}  &  =\sqrt{\frac{\max(\gamma,\alpha_{v})}{\beta_{v}}},\\
L_{c}  &  =\sqrt{\frac{2\gamma}{\min(\alpha_{h}-\alpha_{v},\beta_{v}-\beta
_{h})}}.
\end{align*}
Thus our method gives explicit Lipschitz type bounds for our invariant manifolds of (\ref{eq:our-ode}).
\end{corollary}



\section{Centre manifold around $L_{1}$ in the Restricted Three body problem}

\label{sec:3bp}In the following we specialise our study to the center manifold
of the equilibrium point $L_{1}$ in the restricted three body problem, or RTBP
for short.

Section~\ref{sec:RTBP} describes the RTBP and presents its equations of motion
and specifies the equilibrium point $L_{1}$ around which we shall later prove
existence of the center manifold. A general reference for this section is
Szebehely's book~\cite{Szeb}. Section \ref{sec:NormalForm} constructs
\textquotedblleft aligned coordinates\textquotedblright\ (described in Section
\ref{sec:adjusted-coordinates}) around $L_{1}$ in the RTBP using a suitable
normal form procedure. A general reference for this section is the paper by
Jorba~\cite{Jorba} on computation of normal forms with application to the
RTBP. In Section \ref{sec:CenterManifold} we show how normal forms can be used
to obtain a very accurate numerical estimate on where the centre manifold is
positioned. In Section \ref{sec:application} we apply Theorem
\ref{th:main} to obtain a rigorous enclosure of the centre manifold.

\subsection{Restricted Three Body Problem}

\label{sec:RTBP}

The problem is defined as follows: two main bodies rotate in the plane about
their common center of mass on circular orbits under their mutual
gravitational influence. A third body moves in the same plane of motion as the
two main bodies, attracted by the gravitation of previous two but not
influencing their motion. The problem is to describe the motion of the third body.

Usually, the two rotating bodies are called the \emph{primaries}. We will
consider as primaries the Sun and the Earth. The third body can be regarded as
a satellite or a spaceship of negligible mass.

We use a rotating system of coordinates centred at the center of mass. The
plane $X,Y$ rotates with the primaries. The primaries are on the $X$ axis, the
$Y$ axis is perpendicular to the $X$ axis and contained in the plane of rotation.

\begin{figure}[ptb]
\begin{center}
\includegraphics[height=5cm]{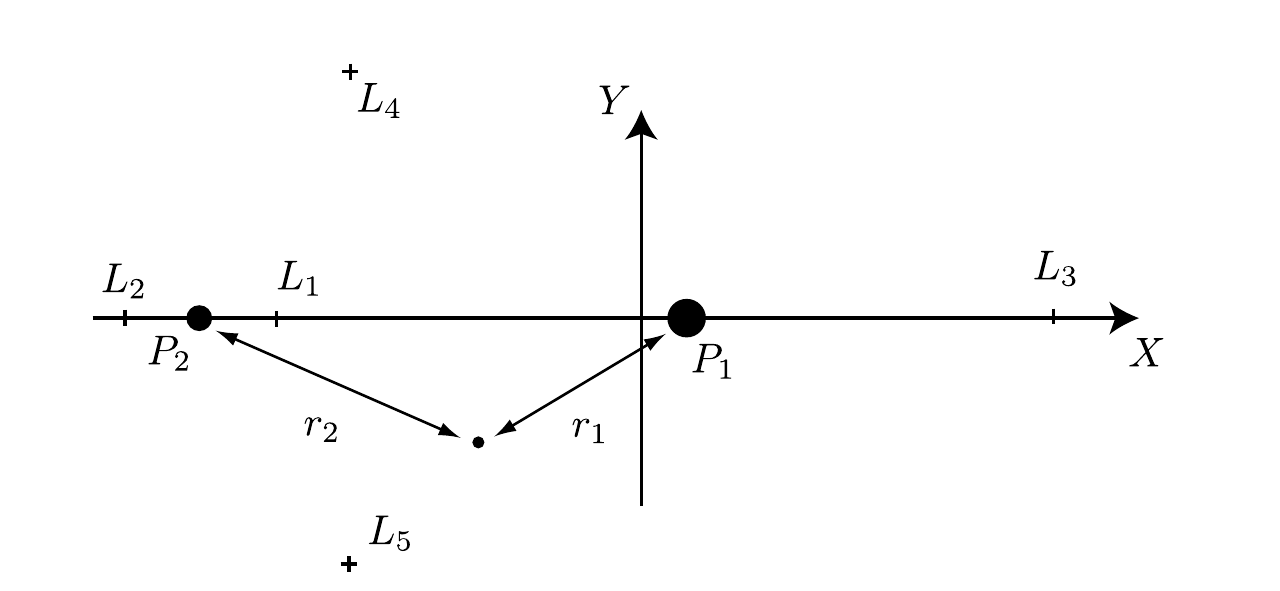}
\end{center}
\caption{Notation for the rotating system of coordinates with origin at the
center of mass. The Sun has the mass $1-\mu$ and is fixed at $P_{1}=(\mu,0)$.
The Earth has the mass $\mu$ is fixed at $P_{2}=(\mu-1,0)$. The third massless
particle moves in the $XY$ plane.}%
\label{fig:3bp}%
\end{figure}

We rescale the masses $\mu_{1}$ and $\mu_{2}$ of the primaries so that they
satisfy the relation $\mu_{1}+\mu_{2}=1$. After such rescaling the distance
between the primaries is $1$. (See Szebehelly~\cite{Szeb}, section 1.5).

Let the smaller mass be $\mu_{2}=\mu=3.040423398444176\times10^{-6}$ and the
larger one be $\mu_{1}=1-\mu$, corresponding to the values of the Earth and
the Sun respectively. We use a convention in which in the rotating coordinates
the Sun is located to the right of the origin at $P_{1}=(\mu,0)$, and the
Earth is located to the left at $P_{2}=(\mu-1,0)$.

The equations of motion of the third body are
\begin{subequations}
\label{eq:RTBP}%
\begin{align*}
\ddot{X}-2\dot{Y}  &  =\Omega_{X},\\
\ddot{Y}+2\dot{X}  &  =\Omega_{Y},
\end{align*}
where
\end{subequations}
\[
\Omega=\frac{1}{2}(X^{2}+Y^{2})+\frac{1-\mu}{r_{1}}+\frac{\mu}{r_{2}},
\]
and $r_{1},r_{2}$ denote the distances from the third body to the larger and
the smaller primary, respectively (see Figure \ref{fig:3bp})
\begin{align*}
r_{1}^{2}  &  =(X-\mu)^{2}+Y^{2},\\
r_{2}^{2}  &  =(X-\mu+1)^{2}+Y^{2}.
\end{align*}
These equations have an integral of motion~\cite{Szeb} called the Jacobi
integral
\[
C=2\Omega-(\dot{X}^{2}+\dot{Y}^{2}).
\]

The equations of motion take Hamiltonian form if we consider positions $X$,
$Y$ and momenta $P_{X}=\dot{X}-Y$, $P_{Y}=\dot{Y}+X$. The Hamiltonian is
\begin{equation}
H=\frac{1}{2}(P_{X}^{2}+P_{Y}^{2})+YP_{X}-XP_{Y}-\frac{1-\mu}{r_{1}}-\frac
{\mu}{r_{2}}, \label{eq:Hamiltonian}%
\end{equation}
with the vector field given by%
\begin{align*}
F  &  =J\nabla H,\\
J  &  =\left(
\begin{array}
[c]{cc}%
0 & \mathrm{id}\\
-\mathrm{id} & 0
\end{array}
\right)  ,\quad\mathrm{id}=\left(
\begin{array}
[c]{cc}%
1 & 0\\
0 & 1
\end{array}
\right)  .
\end{align*}
The Hamiltonian and the Jacobi integral are simply related by $H=-\frac{C}{2}$.

Due to the Hamiltonian integral, the dimensionality of the space can be
reduced by one. Trajectories of equations~\eqref{eq:RTBP} stay on the
\emph{energy surface} $M$ given by $H(X,Y,P_{X},P_{Y})=h=constant$, a
3-dimensional submanifold of $\mathbb{R}^{4}$. Equivalently, $M$ is the level
surface
\begin{equation}
M\equiv\{C(X,Y,\dot{X},\dot{Y})=c=-2h\} \label{eq:energySurface}%
\end{equation}
of the Jacobi integral.

The restricted three body problem in a rotating frame, described by
equations~\eqref{eq:RTBP}, has five equilibrium points (see \cite{Szeb}).
Three of them, denoted $L_{1},L_{2}$ and $L_{3}$, lie on the X axis and are
usually called the `collinear' equilibrium points (see Figure \ref{fig:3bp}).
Notice that we denote $L_{1}$ the \emph{interior} collinear point, located
between the primaries.

At this point we would like to make it clear that in this paper we focus only
on the equilibrium point $L_{1}$, though other collinear equilibria points
could be investigated in the same manner.

The Jacobian of the vector field at $L_{1}$ has two real and two purely
imaginary eigenvalues. Since the three body problem is Hamiltonian in can be
shown by the Lyapunov-Moser theorem \cite{Moser} that in a sufficiently small
neighbourhood of $L_{1}$ there exists a family of periodic orbits which is
parameterised by energy. This family of orbits forms a center manifold. Our
aim shall be to prove the existence of this manifold in a neighbourhood which
is far from $L_{1}.$ As mentioned before, close to $L_{1}$ the existence of
this manifold follows from the center manifold theorem (or in this case also
from the Lyapunov-Moser theorem). The hard task is to prove its existence far
from the equilibrium point.

\begin{remark}\label{rem:LapMoser}
Since the center manifold around $L_{1}$ is foliated by periodic orbits it has
to be identical to the invariant manifold obtained through Theorem
\ref{th:main} due to point 2. of the theorem. The Lyapunov-Moser theorem
ensures the existence of periodic orbits locally. In such local domain we are
guaranteed that the manifold $\chi$ from Theorem \ref{th:main} inherits all
regularity properties which follow from the center manifold theorem. Outside
of this domain Theorem \ref{th:main} establishes only Lipschitz continuity of $\chi$. 
\end{remark}



\subsection{Normal Form}

\label{sec:NormalForm}

The linearised dynamics around the equilibrium point $L_{1}$ is of type saddle
$\times$ centre for all values of $\mu$. In this section we use a normal form
procedure to approximate the nonlinear dynamics locally around $L_{1}$.

For the purpose of this paper, the normal form coordinates will be used
precisely as the well-aligned coordinates described in
section~\ref{sec:adjusted-coordinates}.

The goal of the normal form procedure is to simplify the Taylor expansion of
the Hamiltonian around the equilibrium point using canonical, near-identity
changes of variables. This procedure is carried up to a given (finite) degree
in the expansion. The resulting Hamiltonian is then truncated to (finite)
degree. Such Hamiltonian is said to be in~\emph{normal form}.

We compute a normal form expansion that is as simple as possible, i.e. one
that has the minimum number of monomials. This is sometimes called a
\emph{full}, or \emph{complete}, normal form. The equations of motion
corresponding to the truncated normal form can be integrated exactly. As a
result, \emph{locally} the normal form gives a very accurate approximation of
the dynamics.

In particular, here we use the normal form to approximate the local center
manifold by a 1-parameter family of periodic orbits with increasing energy.

The normal form construction proceeds in three steps. First we perform some
convenient translation and scaling of coordinates, and expand the Hamiltonian
around $L_{1}$ as a power series. Then we make a linear change of coordinates
to put the quadratic part of the Hamiltonian in a simple form, which
diagonalises the linear part of equations of motion. Finally we use the
so-called Lie series method to perform a sequence of canonical, near-identity
transformations that simplify nonlinear terms in the Hamiltonian of
successively higher degree.

The transformation to well-aligned coordinates $\phi:U\rightarrow
\phi(U)\subset\mathbb{R}^{n}$ is the composition of all the transformations
performed during these three steps.

A similar full normal form expansion has been used for the \emph{spatial} RTBP
in a previous paper~\cite{DelshamsMR08}. We refer to the previous paper for
the fine details of the normal form construction, which will be left out of
the current paper.

\subsubsection{Hamiltonian expansion}

We start by writing the Hamiltonian~\eqref{eq:Hamiltonian} as a power series
expansion around the equilibrium point $L_{1}$. First we translate the origin
of coordinates to the equilibrium point. In order to have good numerical
properties for the Taylor coefficients, it is also convenient to scale
coordinates~\cite{Richardson80}. The translation and scaling are given by
\begin{equation}
X=-\gamma x+\mu-1+\gamma,\quad Y=-\gamma y, \label{eq:scaling}%
\end{equation}
where $\gamma$ is the distance from $L_{1}$ to its closest primary (the Earth).

Since scalings are not canonical transformations, we apply this change of
coordinates to the equations of motion, to obtain
\begin{subequations}
\label{eq:scaledRTBP}%
\begin{align*}
\ddot{x}-2\dot{y}  &  =\Omega_{x}\\
\ddot{y}+2\dot{x}  &  =\Omega_{y},
\end{align*}
where
\end{subequations}
\[
\Omega=\frac{1}{2}(x^{2}+y^{2})-\frac{\mu-1+\gamma}{\gamma}x+\frac{1}%
{\gamma^{3}}\left(  \frac{1-\mu}{r_{1}}+\frac{\mu}{r_{2}}\right)
\]
and $r_{1},r_{2}$ denote the (scaled) distances from the third body to the
larger and the smaller primary, respectively.

Defining $p_{x}=\dot{x}-y$, $p_{y}=\dot{y}+x$, the libration-point centred
equations of motion~\eqref{eq:scaledRTBP} are Hamiltonian, with Hamiltonian
function
\begin{equation}
H=\frac{1}{2}(p_{x}^{2}+p_{y}^{2})+yp_{x}-xp_{y}+\frac{\mu-1+\gamma}{\gamma
}x-\frac{1}{\gamma^{3}}\left(  \frac{1-\mu}{r_{1}}+\frac{\mu}{r_{2}}\right)  .
\label{eq:scaledHamiltonian}%
\end{equation}
Our first change of coordinates can therefore be summarised as $R:\mathbb{R}%
^{4}\rightarrow\mathbb{R}^{4}$%
\begin{align}
\left(  X,Y,P_{X},P_{Y}\right)   &  =R(x,y,p_{x},p_{y}%
)\label{eq:rescaling-change}\\
&  =\left(  -\gamma x+\mu-1+\gamma,-\gamma y,-\gamma p_{x},-\gamma p_{y}%
+\mu-1+\gamma\right) \nonumber
\end{align}

The Hamiltonian is then rewritten in the form \cite{Jorba,JorbaM99}
\begin{equation}
H=\frac{1}{2}(p_{x}^{2}+p_{y}^{2})+yp_{x}-xp_{y}-\sum_{n\geq2}c_{n}(\mu
)\rho^{n}P_{n}\left(  \frac{x}{\rho}\right)  , \label{eq:HamiltonianExpansion}%
\end{equation}
where $P_{n}$ is the $n$-th Legendre polynomial, and the coefficients
$c_{n}(\mu)$ are given by
\[
c_{n}(\mu)=\frac{1}{\gamma^{3}}\left(  \mu+(-1)^{n}\frac{(1-\mu)\gamma^{n+1}%
}{(1-\gamma)^{n+1}}\right)  .
\]
This expansion holds when $\rho<\min{(|P_{1}|,|P_{2}|)}=|P_{2}|=1$, i.e. it
is valid in a ball centred at $L_{1}$ that extends up to the Earth.

\subsubsection{Linear changes of coordinates}

Now we transform the \emph{linear} part of the system into Jordan form, which
is convenient for the normal form procedure. This particular transformation is
derived in \cite{Jorba,JorbaM99}, for instance.

Consider the quadratic part $H_{2}$ of the
Hamiltonian~\eqref{eq:HamiltonianExpansion},
\begin{equation}
H_{2}=\frac{1}{2}(p_{x}^{2}+p_{y}^{2})+yp_{x}-xp_{y}-c_{2}x^{2}+\frac{c_{2}%
}{2}y^{2}, \label{eq:H2}%
\end{equation}
which corresponds to the linearisation of the equations of motion. It is
well-known~\cite{JorbaM99} that the linearised system has eigenvalues of the
form $\pm\lambda,\pm i\nu$, where $\lambda,v$ are real and positive.

One can find (\cite{JorbaM99} section~2.1) a symplectic linear change of
variables
\[
C=\left(
\begin{array}
[c]{cccc}%
\frac{2\lambda}{s_{1}} & \frac{-2\lambda}{s_{1}} & 0 & \frac{2v}{s_{2}}\\
\frac{\lambda^{2}-2c_{2}-1}{s_{1}} & \frac{\lambda^{2}-2c_{2}-1}{s_{1}} &
\frac{-v^{2}-2c_{2}-1}{s_{2}} & 0\\
\frac{\lambda^{2}+2c_{2}+1}{s_{1}} & \frac{\lambda^{2}+2c_{2}+1}{s_{1}} &
\frac{-v^{2}+2c_{2}+1}{s_{2}} & 0\\
\frac{\lambda^{3}+(1-2c_{2})\lambda}{s_{1}} & \frac{-\lambda^{3}%
-(1-2c_{2})\lambda}{s_{1}} & 0 & \frac{-v^{3}+(1-2c_{2})v}{s_{2}}%
\end{array}
\right)
\]
where%
\begin{align*}
s_{1}  &  =\sqrt{2\lambda\left(  \left(  4+3c_{2}\right)  \lambda^{2}%
+4+5c_{2}-6c_{2}^{2}\right)  },\\
s_{2}  &  =\sqrt{v\left(  \left(  4+3c_{2}\right)  v^{2}-4-5c_{2}+6c_{2}%
^{2}\right)  },
\end{align*}
that puts the linear terms of the vector field at $L_{1}$ into a Jordan form.
This means that the change from position-momenta to new variables
$(x_{1},y_{1},x_{2},y_{2})\in\reals^{4}$,
\begin{equation}
(x,y,p_{x},p_{y})=C(x_{1},y_{1},x_{2},y_{2}), \label{eq:linearChange}%
\end{equation}
casts the quadratic part of the Hamiltonian into
\begin{equation}
H_{2}=\lambda x_{1}y_{1}+\frac{\nu}{2}(x_{2}^{2}+y_{2}^{2}). \label{eq:H2_nf}%
\end{equation}

The linear equations of motion $(\dot{x},\dot{y})=A(x,y)$ associated
to~\eqref{eq:H2_nf} decouple into a hyperbolic and a center part
\begin{subequations}
\label{eq:linearNF}%
\begin{align}
(\dot{x}_{1},\dot{y}_{1}) &  =A_{h}(x_{1},y_{1})\\
(\dot{x}_{2},\dot{y}_{2}) &  =A_{c}(x_{2},y_{2}),
\end{align}
with
\end{subequations}
\[
A_{h}=%
\begin{pmatrix}
\lambda & 0\\
0 & -\lambda
\end{pmatrix}
\qquad\text{and }\qquad A_{c}=\left(
\begin{array}
[c]{cc}%
0 & v\\
-v & 0
\end{array}
\right)  .
\]

Notice that the matrix $A$ of the linear equations (\ref{eq:linearNF}) is in
block-diagonal form. It is convenient to diagonalise the matrix $A$ over
$\mathbb{C.}$ Consider the \emph{symplectic} change $T^{-1}\colon
\reals^{4}\rightarrow\complexs^{4}$ to complex variables $(q,p)=(q_{1}%
,p_{1},q_{2},p_{2})\in\complexs^{4}$%
\begin{align}
(q_{1},p_{1},q_{2},p_{2})  &  =T^{-1}(x_{1},y_{1},x_{2},y_{2}%
)\label{eq:to-complex-change}\\
&  =\left(  x_{1},y_{1},\frac{1}{\sqrt{2}}(x_{2}-iy_{2}),\frac{1}{\sqrt{2}%
}(-ix_{2}+y_{2})\right)  .\nonumber
\end{align}
This change casts the quadratic part of the Hamiltonian into
\begin{equation}
H_{2}=\lambda q_{1}p_{1}+i\nu q_{2}p_{2}. \label{eq:H2_complex_nf}%
\end{equation}
Equivalently, this change carries $A$ to diagonal form:
\[
T^{-1}AT=\Lambda=\mathrm{diag}(\lambda,-\lambda,i\nu,-i\nu).
\]

\subsubsection{Nonlinear normal form}

Assume that the symplectic linear changes of variables~\eqref{eq:linearChange}
and (\ref{eq:to-complex-change}) have been performed in the Hamiltonian
expansion~\eqref{eq:HamiltonianExpansion}, so that the quadratic part $H_{2}$
is already in the form~\eqref{eq:H2_complex_nf}.

Let us thus write the Hamiltonian as
\begin{equation}
H(q,p)=H_{2}(q,p)+H_{3}(q,p)+H_{4}(q,p)+\cdots\label{eq:HamiltonianSeries}%
\end{equation}
with $H_{j}(q,p)$ as homogeneous polynomials of degree $j$ in the variables
$(q,p)\in\complexs^{4}$.

As shown in a previous paper~\cite{DelshamsMR08}, we can remove most monomials
in the series~\eqref{eq:HamiltonianSeries} by means of formal coordinate
transformations, in order to obtain an integrable approximation to the
dynamics close to the equilibrium point.

\begin{proposition}
[Complete normal form around a saddle$\times$centre]%
\label{prop:normal_form}

\cite{DelshamsMR08} For any integer $N\geq3$, there exists a neighbourhood $\mathcal{U}^{(N)}$ of
the origin and a near-identity canonical transformation
\begin{equation}
\mathcal{T}^{(N)}\colon\ \complexs^{4}\supset\mathcal{U}^{(N)}\mapsto
\complexs^{4} \label{eq:nonlinear-nf}%
\end{equation}
that puts the system~\eqref{eq:HamiltonianSeries} in normal form up to order
$N$, namely
\[
\mathcal{H}^{(N)}:=H\circ\mathcal{T}^{(N)}=H_{2}+\mathcal{Z}^{(N)}%
+\mathcal{R}^{(N)}%
\]
where $\mathcal{Z}^{(N)}$ is a polynomial of degree $N$ that Poisson-commutes
with $H_{2}$
\[
\{\mathcal{Z}^{(N)},H_{2}\}\equiv0,
\]
and $\mathcal{R}^{(N)}$ is small
\[
|\mathcal{R}^{(N)}(z)|\leq C_{N}||z||^{N+1}\quad\forall z\in\mathcal{U}%
^{(N)}.
\]

If the elliptic frequencies $\nu,\omega$ are nonresonant to degree $N$,
\[
c_{1}\nu+c_{2}\omega\neq0\qquad\forall(c_{1},c_{2})\in\integers^{2}%
,\quad0<|c_{1}+c_{2}|\leq N,
\]
then in the new coordinates, the truncated Hamiltonian $H_{2}+\mathcal{Z}%
^{(N)}$ depends only on the \emph{basic invariants}
\begin{subequations}
\begin{align}
I_{1} &  =q_{1}p_{1}=x_{1}y_{1}\label{eq:hyperbolic_coords}\\
I_{2} &  =iq_{2}p_{2}=q_{2}\bar{q}_{2}=\frac{x_{2}^{2}+y_{2}^{2}}%
{2}.\label{eq:elliptic_coords}%
\end{align}

\end{subequations}
\end{proposition}

The equations of motion associated to the \emph{truncated normal form}
$H_{2}+\mathcal{Z}^{(N)}$ can be integrated exactly.

\begin{remark}
The reminder $\mathcal{R}^{(N)}$ is very small in a small neighbourhood of the
origin. Hence, close to the origin, the exact solution of the truncated normal
form is a very accurate approximate solution of the original system $H$.
\end{remark}

\begin{remark}
Let $\phi_{1},\phi_{2}$ be the symplectic conjugate variables to $I_{1},I_{2}%
$, respectively. The basic invariant $I_{2}$ is usually called \emph{action
variable}, and its conjugate variable $\phi_{2}$ is usually called \emph{angle
variable}. They are given in polar variables (\ref{eq:elliptic_coords}).
\end{remark}

We can now write our function $\phi$  for our
change into the well \emph{aligned coordinates} (\ref{eq:phi-total-change}). To do so we compose the inverse
transformations given in (\ref{eq:rescaling-change}), (\ref{eq:linearChange}),
(\ref{eq:to-complex-change}) and (\ref{eq:nonlinear-nf}) which gives us%
\begin{equation}
\phi=\left(  \mathcal{T}^{(N)}\right)  ^{-1}\circ T^{-1}\circ C^{-1}\circ
R^{-1}. \label{eq:nf-phi}%
\end{equation}

\begin{remark}
The above described method of obtaining normal form coordinates is performed
by passing through complex variables. It is possible though to arrange the
changes so that the combined change of coordinates (\ref{eq:nf-phi}) passes
from real to real coordinates. The change of coordinates $\phi$ is a high
order polynomial. It is possible to arrange the normal form change of
coordinates so that the coefficient of $\phi$ are real (see \cite{Jorba}). In
setting up our change of coordinates for the application of Theorem
\ref{th:main} to the RTBP in Section
\ref{sec:application} we have adopted such a procedure.
\end{remark}

\begin{remark}
In practice, one usually computes a normal form of degree $N=16$. In our
application to the restricted three body problem in Section
\ref{sec:application} we use a normal form of degree $N=4$. This turns out to
be sufficient, since we investigate a relatively close neighbourhood of the
invariant point, where degree of order four gives us a sufficiently good approximation.
\end{remark}



\subsection{Approximating the center manifold in normal form coordinates}

\label{sec:CenterManifold}

In normal form coordinates given by (\ref{eq:nf-phi}) the Hamiltonian, by
Proposition \ref{prop:normal_form}, is of the form%
\begin{equation}
\mathcal{H}^{(N)}=H_{2}+\mathcal{Z}^{(N)}+\mathcal{R}^{(N)},\qquad
\{\mathcal{Z}^{(N)},H_{2}\}\equiv0.\label{eq:NF-order-N}%
\end{equation}
In this section we shall show that when we neglect the reminder term
$\mathcal{R}^{(N)}$, and thus consider an approximation of the system, the
normal form coordinates given by (\ref{eq:nf-phi}) give us a very good
understanding of where the center manifold is positioned and of the dynamics
on it.

Let $U$ be some small neighbourhood of the fixed point (in our discussion for
the R3BP this will be $L_{1}$) and let $\phi\colon U\rightarrow\phi
(U)\subset\reals^{4}$ be the transformation to normal form coordinates
(\ref{eq:nf-phi}). Consider the normal form (\ref{eq:NF-order-N}) up to order
$N$ with associated equations of motion
\begin{equation}
\dot{p}=F^{\phi}(p):=J\nabla\mathcal{H}^{(N)}(p). \label{eq:full_system}%
\end{equation}
Consider now the \emph{truncated} normal form up to order $N$
\[
\hat{\mathcal{H}}^{(N)}=H_{2}+\mathcal{Z}^{(N)},
\]
with associated equations of motion
\begin{equation}
\dot{p}=\hat{F}^{\phi}(p):=J\nabla\hat{\mathcal{H}}^{(N)}(p).
\label{eq:truncated_system}%
\end{equation}

Recall that the corresponding linearisation around the origin is
(\ref{eq:linearNF})
\begin{equation}
\dot{p}=Ap,\qquad p=(x_{1},y_{1},x_{2},y_{2})\in\reals^{4}%
,\label{eq:linearized_system}%
\end{equation}
where $x_{1},y_{1}$ are the hyperbolic normal form
coordinates~\eqref{eq:hyperbolic_coords}, and $x_{2},y_{2}$ are the center
normal form coordinates~\eqref{eq:elliptic_coords}. In order to match the
notation from section~\ref{sec:setup}, let us denote the center normal form
coordinates $x_{2},y_{2}$ as $\theta_{1},\theta_{2}$, the unstable normal form
coordinate $x_{1}$ as $x$, and the stable normal form coordinate $y_{1}$ as
$y$. Note that to match the notations we need to swap the order in which the
coordinates are written out passing from $(x_{1},y_{1},x_{2},y_{2})$ to
$(\theta_{1},\theta_{2},x,y)$.

The truncated system $\hat{F}^{\phi}$ has several invariant subspaces.
Specifically, the next proposition follows from~\cite{Murdock}, section 5.1.

\begin{proposition}
Let
\begin{align}
E^{c}  &  =\{(\theta_{1},\theta_{2},0,0)\colon\ (\theta_{1},\theta_{2}%
)\in\reals^{2}\},\\
E^{u}  &  =\{(0,0,x,0)\colon\ x\in\reals\},\\
E^{s}  &  =\{(0,0,0,y)\colon\ y\in\reals\}.
\end{align}
Then, $E^{c}$, $E^{u}$ and $E^{s}$ are invariant subspaces of the flow of
$\hat{F}^{\phi}$.
\end{proposition}

\begin{remark}
These subspaces are invariant of the nonlinear \textbf{truncated }system
(\ref{eq:truncated_system}), not just of the linearised
system~\eqref{eq:linearized_system}. It is important to stress here though
that these subspaces need not  be invariant under the full system
(\ref{eq:full_system}).
\end{remark}

Next we claim that $E^{c}$ is approximately equal to the center manifold
$W^{c}$ of the \emph{full system} $F^{\phi}$. This is formulated in the next
proposition which follows from~\cite{Murdock}, Section 5.2.

\begin{proposition}
\label{prop:cmfld-nf}For each integer $r$ with $N\leq r<\infty$, there exists
a (not necessarily unique) local invariant center manifold $W^{c}$ of
$F^{\phi}$ of class $C^{r}$ such that

\begin{itemize}
\item $W^{c}$ is expressible as a graph over $E^{c}$, i.e. there exists a
neighbourhood $V\subset E^{c}$ and a map $\chi\colon V\to E^{u}\oplus E^{s}$
such that
\[
W^{c} = \{ (\theta_{1},\theta_{2},x,y)\in E^{c} \oplus E^{u} \oplus
E^{s}\colon\ (\theta_{1},\theta_{2})\in V, \ (x,y)=\chi(\theta_{1},\theta
_{2})\}.
\]

\item $W^{c}$ has $N$-th order contact with $E^{c}$, i.e. $\chi$ and its
derivatives up to order $N$ vanish at the origin.
\end{itemize}
\end{proposition}

Hence, in normal form coordinates, the center manifold $W^{c}$ of $F^{\phi}$
is approximated very accurately (to order $N$) around the origin by the
subspace $E^{c}$.

\begin{remark}
When applying Proposition \ref{prop:cmfld-nf} we are faced with a problem that
it is usually very hard to obtain a rigorous bound on the size of the set $V$.
Moreover, even though we know that up to order $N$ the derivatives of $\chi$
vanish at zero, it is usually very hard to obtain a rigorous bound on the size
of the higher order terms of $\chi$ on the set $V,$ and thus obtain a rigorous
bound on the position of the true centre manifold.
\end{remark}

Let us now briefly discuss the dynamics of the system
(\ref{eq:truncated_system}) on $E^{c}.$ To do so we shall use the center
normal form coordinates~\eqref{eq:elliptic_coords} in action-angle form, i.e.
from now on we will use $(I,\varphi)\in\reals\times\tori$ for the center part.
Proposition~\ref{prop:normal_form} states that the truncated Hamiltonian
$\hat{\mathcal{H}}^{(N)}$ depends only on the action $I$, and not on the angle
$\varphi$. Thus the restriction of $\hat{F}^{\phi}$ to its invariant subspace
$E^{c}$ is
\begin{equation}
\dot{I}=0,\qquad\dot{\varphi}=\frac{\partial\hat{\mathcal{H}}^{(N)}}{\partial
I}=:\omega(I).\label{eq:dyn-action-angle}%
\end{equation}
The solutions inside $E^{c}$ with initial condition $I(0)=I_{0}$ and
$\varphi(0)=\varphi_{0}$ are $I(t)=I_{0}$, $\varphi(t)=\omega(I_{0}%
)t+\varphi_{0}$. In the case of the restricted three body problem $E^{c}$ is
two dimensional and so the dynamics of the truncated system on $E^{c}$ is
foliated by invariant circles of increasing action $I$. Notice from
equation~\eqref{eq:H2_nf} that $H$ grows linearly with respect to $I$ (close
to the origin), so the invariant circles also have increasing energy $H$.

The properties discussed above motivate the use of the normal form coordinates
$\theta_{1},\theta_{2},x,y$ as the well-aligned coordinates in the sense of
section~\ref{sec:adjusted-coordinates}. They provide a good guess on the
location of the center manifold (locally around the origin). Taking $\bar
{B}_{c}^{R}=\{(\theta_{1},\theta_{2})\in\reals^{2}\colon\ \left\Vert
\theta_{1},\theta_{2}\right\Vert \leq R\}$ the guess is given by $\phi
^{-1}(\bar{B}_{c}^{R}\times\{0\})\subset\reals^{4}$. Notice also that the
center coordinate $I$ is well-aligned with the energy (in sense of (\ref{eq:h-boundary-bound})). Let
\[
C^{R}=\left\{  (\theta_{1},\theta_{2})\in\reals^{2}\colon\ \left\Vert
\theta_{1},\theta_{2}\right\Vert =R\right\}  =\left\{  \left(  I,\varphi
\right)  \in\reals\times\tori:I=\frac{R^{2}}{2}\right\}
\]
be the invariant circle of radius $R$ for the system
(\ref{eq:truncated_system}).\textbf{ }By (\ref{eq:dyn-action-angle}) we have
$\hat{\mathcal{H}}^{(N)}(C^{R_{1}}\times\{0\})<\hat{\mathcal{H}}%
^{(N)}(C^{R_{2}}\times\{0\})$ whenever $R_{1}<R_{2}$. Hence, given an energy $h$, we can find $R_{1},R_{2}>0$ such that
\[
\hat{\mathcal{H}}^{(N)}(C^{R_{1}}\times\{0\})<h<\hat{\mathcal{H}}%
^{(N)}(C^{R_{2}}\times\{0\}).
\]
Taking $R_{1},R_{2}$ sufficiently far (in practice they are still close) from
one another and taking sufficiently small $r>0$, since $\hat{\mathcal{H}%
}^{(N)}$ and $\mathcal{H}^{(N)}$ are close, we expect also that
\[
\mathcal{H}^{(N)}(B_{c}^{R_{1}}\times B_{u}^{r}\times B_{s}^{r})<h<\mathcal{H}%
^{(N)}(C^{R_{2}}\times B_{u}^{r}\times B_{s}^{r}).
\]
Since $\mathcal{H}^{(N)}=H\circ\phi^{-1}$, this will mean that that the
bound~\eqref{eq:h-boundary-bound} shall be satisfied.



\subsection{Application of the main theorem to the center manifold around
$L_{1}$}

\label{sec:application}

In this section we shall show how to apply Theorem \ref{th:main} in practice.

As described in Section \ref{sec:NormalForm} the change of coordinates to
\emph{well aligned coordinates }can be done using a change to normal
coordinates (\ref{eq:nf-phi}). We obtain the function $\phi$ using the
algorithm of Jorba \cite{Jorba}. The algorithm allows us to obtain
$\phi$ as a real polynomial, passing from $\mathbb{R}^{4}$ to $\mathbb{R}^{4}$.

\subsubsection{Methodology}

To apply Theorem \ref{th:main} it is enough to derive a rigorous bound on the
derivative of $F^{\phi}$. Let us now outline how such a bound can be obtained.
Using (\ref{eq:field-local}), for any $p\in\mathbb{R}^{4}$ we have%
\begin{align}
D\left(  F^{\phi}(p)\right)   &  =D^{2}\phi(\phi^{-1}(p))D(\phi^{-1}%
)(p)F(\phi^{-1}(p))\label{eq:DFphi-at-p}\\
&  +D\phi(\phi^{-1}(p))DF(\phi^{-1}(p))D(\phi^{-1})(p).\nonumber
\end{align}
In our computer assisted proof we apply the above formula using an
interval-arithmetic-based software called
CAPD (Computer Assisted Proofs in
Dynamics\footnote{http://capd.ii.uj.edu.pl}). This software in particular allows for
rigorous-interval-enclosure-based computation of high order derivatives of
functions on sets. In our application we obtain a global bound for the derivative (\ref{eq:dF-bound}) on the entire set $D_{\phi}$. To
compute $[DF^{\phi}(D_{\phi})]$ applying (\ref{eq:DFphi-at-p}) we only requite
to compute images of functions, derivatives of functions and a second
derivative on a set $D_{\phi}$. All such computations can be performed in CAPD.

Before specifying the size of the set $D_{\phi}$ and giving
rigorous-interval-based numerical results, we have to stress one technical
problem encountered when applying formula (\ref{eq:DFphi-at-p}). We take our
change to \emph{well aligned coordinates} $\phi$ to be a high order polynomial
obtained from non-rigorous computations. To apply formula (\ref{eq:DFphi-at-p}%
) directly we would need to know its inverse $\phi^{-1}$. Let us stress that
one can not use a numerical approximation of an inverse change and use it as
$\phi^{-1}$ (such numerical approximate inverse is readily available from
algorithms of \cite{Jorba}). To apply (\ref{eq:DFphi-at-p}) directly one would
have to use a \emph{rigorous}, analytic inverse. Since $\phi$ is a polynomial in high dimension and of high order, its analytic inverse is next to impossible to obtain in practice. To remedy this problem we
slightly modify (\ref{eq:DFphi-at-p}). Using the fact that $D(\phi
^{-1})(p)=\left(  D\phi(\phi^{-1}(p))\right)  ^{-1}$ we can rewrite
(\ref{eq:DFphi-at-p}) as%
\begin{align*}
DF^{\phi}(p)  &  =D^{2}\phi(\phi^{-1}(p))(D\phi(\phi^{-1}(p)))^{-1}F(\phi
^{-1}(p))\\
&  +D\phi(\phi^{-1}(p))DF(\phi^{-1}(p))\left(  D\phi(\phi^{-1}(p))\right)
^{-1}.
\end{align*}
This in interval arithmetic notation gives us the following formula for the
interval enclosure of $DF^{\phi}$ on some set $\mathbf{I}\subset D_{\phi}$%
\begin{align}
\lbrack DF^{\phi}(\mathbf{I})]  &  \subset\lbrack D^{2}\phi(\left[  \phi
^{-1}(\mathbf{I})\right]  )(D\phi(\left[  \phi^{-1}(\mathbf{I})\right]
))^{-1}F(\left[  \phi^{-1}(\mathbf{I})\right]  )\label{eq:DF-phi-I}\\
&  +D\phi(\left[  \phi^{-1}(\mathbf{I})\right]  )DF(\left[  \phi
^{-1}(\mathbf{I})\right]  )\left(  D\phi(\left[  \phi^{-1}(\mathbf{I})\right]
)\right)  ^{-1}]\nonumber
\end{align}
To compute the right hand side of the above equation there is no need to
invert the function $\phi.$ It is enough to find a set $\left[  \phi
^{-1}(\mathbf{I})\right]  $ which contains the pre-image of $\mathbf{I}$,
i.e.
\[
\phi^{-1}(\mathbf{I})\subset\left[  \phi^{-1}(\mathbf{I})\right]  ,
\]
and for this we do not need to compute the inverse function. For a set
$B\subset\mathbb{R}^{4}$ the following lemma can be used to verify that
$\phi^{-1}(\mathbf{I})\subset B.$

\begin{lemma}
\label{lem:preimage-phi}Let $\phi:\mathbb{R}^{n}\rightarrow\mathbb{R}^{n}$ be
a homeomorphism and let $\mathbf{I},B\subset\mathbb{R}^{n}$ be two sets
homeomorphic to $n$-dimensional balls. If $\phi(\partial B)\cap\mathbf{I}%
=\emptyset$ and for some point $p\in B$ we have $\phi(p)\in\mathbf{I}$ then%
\[
\phi^{-1}(\mathbf{I})\subset B.
\]

\end{lemma}

\begin{proof}
This follows from elementary topological arguments.
\end{proof}

To apply the lemma in practice it is convenient to first have a non-rigorous guess
on the inverse function, let us denote it by $\hat{\phi}^{-1}$. This means
that
\[
\hat{\phi}^{-1}\phi\approx\mathrm{id.}%
\]
A function $\hat{\phi}^{-1}$ is readily available from algorithms of Jorba
\cite{Jorba}. We can then choose $\lambda>1$ and set $B=\lambda\lbrack
\hat{\phi}^{-1}(\mathbf{I})]$ (in our application we choose $\lambda=3$ which
we find is large enough for our problem). Then we divide the boundary
$\partial B$ into smaller sets and verify that the image by $\phi$ of each
smaller set is disconnected with $\mathbf{I}$. We also check that for the
middle point $p$ in $B$ we have $\phi(p)\in\mathbf{I}$. This by Lemma
\ref{lem:preimage-phi} guarantees that $\phi^{-1}(\mathbf{I})\subset B$.

\begin{remark}
\label{rem:preimage}Once a set $B$ such that $\phi^{-1}(\mathbf{I})\subset B$
is found, there is a useful trick that can be used to refine this initial
guess on the pre-image. One can take a very small set $\mathbf{I}_{0}%
\subset\mathbf{I}$ and using Lemma \ref{lem:preimage-phi} find a small set
$B_{0}$ such that $\phi^{-1}(\mathbf{I}_{0})\subset B_{0}.$ The set $\left[
\phi^{-1}(\mathbf{I})\right]  $ can then be chosen as%
\[
\left[  \phi^{-1}(\mathbf{I})\right]  =B_{0}+\left[  (D\phi(B))^{-1}\right]
[\mathbf{I-I}_{0}].
\]
Such choice guarantees that $\phi^{-1}(\mathbf{I})\subset\left[  \phi
^{-1}(\mathbf{I})\right]  .$ It is also usually tighter than the
initial guess $B$; which is true especially when the function $\phi$ is close
to identity.
\end{remark}

\begin{proof}
This follows from the mean value theorem.
\end{proof}

In a similar fashion to the method from Remark \ref{rem:preimage}, to compute
the energy for a set $\mathbf{I}\subset D^{\phi}$, we take some small set
$\mathbf{I}_{0}\subset\mathbf{I}$ and compute%
\begin{equation}
\lbrack H(\phi^{-1}(\mathbf{I}))]\subset H([\phi^{-1}(\mathbf{I}%
_{0})])+\left[  DH([\phi^{-1}(\mathbf{I})])\right]  \left[  [\phi
^{-1}(\mathbf{I})]-[\phi^{-1}(\mathbf{I}_{0})]\right]  .
\label{eq:Energy-computation}%
\end{equation}

\begin{remark}
When applying above tools to compute $[DF^{\phi}(\mathbf{I})]$ using
(\ref{eq:DF-phi-I}), it pays off to use the fact that $\phi$ is composed of
linear changes of coordinates, together with a nonlinear change $\mathcal{T}%
^{(N)}$ which is close to identity. Keeping track of both linear and nonlinear
change allows to tighten the interval bounds of computations.
\end{remark}

To prove the existence of a fixed point (in case of the RTBP we take the point
$L_{1}$) inside of our set $D_{\phi}$ we use the interval Newton method.

\begin{theorem}
\label{th:interval-Newton}\cite{Alefeld}Let $F:\mathbb{R}^{n}\rightarrow
\mathbb{R}^{n}$ be a $C^{1}$ function. Let $\mathbf{I}=\Pi_{i=1}^{n}%
[a_{i},b_{i}],$ $a_{i}<b_{i}.$ Assume that the interval enclosure of
$DF(\mathbf{I})$, denoted by $\left[  DF(\mathbf{I})\right]  $, is invertible.
Let $x_{0}\in\mathbf{I}$ and define
\[
N(F,x_{0},\mathbf{I})=-\left[  DF\left(  \mathbf{I}\right)  \right]
^{-1}F(x_{0})+x_{0}.
\]
If $N(x_{0},\mathbf{I})\subset\mathbf{I}$ then there exists a unique point
$x^{\ast}\in\mathbf{I}$ such that $F(x^{\ast})=0.$
\end{theorem}

\subsubsection{Rigorous-interval-based numerical results}

For our proof we use a normal form of order $N=4$ change of coordinates (\ref{eq:nf-phi}). At this point we stress once again that $\phi$ obtained
by (\ref{eq:nf-phi}) does not need to perfectly align coordinates. A
numerically obtained polynomial, provided that it aligns coordinates well
enough, is sufficient to prove the existence of a center manifold using our
method, provided that assumptions of Theorem \ref{th:main} can be verified.

We investigate a set
\[
D_{\phi}=\bar{B}_{c}^{R}\times\bar{B}_{u}^{r}\times\bar{B}_{s}^{r}%
\]
with%
\begin{align}
R  &  =\sqrt{2\cdot155\cdot10^{-4}}\approx0.176,\label{eq:R-numerics}\\
r  &  =5\cdot10^{-4}.\nonumber
\end{align}
$\,$Our choice of $R$ by (\ref{eq:hyperbolic_coords}) implies that we consider
actions $I\in\lbrack0,0.0155].$

We first prove that we have a fixed point in $D^{\phi}$ applying Theorem
\ref{th:interval-Newton}. We take $\mathbf{I}=\Pi_{i=1}^{4}[-25\cdot
10^{-5},25\cdot10^{-5}]\subset\mathrm{int}D^{\phi}$ and compute
\begin{multline*}
N(F^{\phi},0,\mathbf{I}) \subset \\
\subset \mathtt{[-6.24567e-05,6.245664e-05]}\times\mathtt{[-6.24434e-05,6.24435e-05]}%
\\
\times\mathtt{[-6.24908e-05,6.24908e-05]}\times
\mathtt{[-5.33554e-05,5.33554e-05]}.
\end{multline*}
Clearly $N(F^{\phi},0,\mathbf{I})\subset\mathbf{I}$, which establishes that $L_{1}$ is in the interior of $D^{\phi}$.

Next we verify condition (\ref{eq:h-boundary-bound}). We take $v=\sqrt
{2\cdot155\cdot10^{-4}}-\sqrt{2\cdot150\cdot10^{-4}},$ which is equivalent to
the ball $B_{c}^{R-v}$ having actions $I\in\lbrack0,0.015]$. We subdivide
$\bar{B}_{u}^{r}\times\bar{B}_{s}^{r}$ into $9$ pieces and cover $\partial
\bar{B}_{c}^{R}$ by $500$ small boxes in $\mathbb{R}^2$. Taking the $9\cdot500$ sets, using
(\ref{eq:Energy-computation}) we obtain a bound on the energy
\[
H(\phi^{-1}(\partial\bar{B}_{c}^{R}\times\bar{B}_{u}^{r}\times\bar{B}_{s}%
^{r}))>\mathtt{-1.500445781623899.}%
\]
Taking same type of subdivisions we then show that%
\[
H(\phi^{-1}(\bar{B}_{c}^{R-v}\times\bar{B}_{u}^{r}\times\bar{B}_{s}%
^{r}))<\mathtt{-1.500445787982231,}%
\]
which establishes (\ref{eq:h-boundary-bound}).

To compute $[DF^{\phi}(D_{\phi})],$ we cover $\bar{B}_{c}^{R}$ by $31\,000$
boxes in $\mathbb{R}^{2}$ , $\bar{B}_{c}^{R}\subset%
{\textstyle\bigcup_{i=1}^{31\,000}}
\mathbf{I}_{c,i}.$ Taking $\mathbf{I}_{i}=\mathbf{I}_{c,i}\times\bar{B}%
_{u}^{r}\times\bar{B}_{s}^{r}$ we compute the bound on $[DF^{\phi}(D_{\phi})]$
as an interval hull of all matrices $[DF^{\phi}(\mathbf{I}_{i})]$ (This means
that we take an interval matrix $[DF^{\phi}(D_{\phi})]$ so that $[DF^{\phi
}(\mathbf{I}_{i})]\subset\lbrack DF^{\phi}(D_{\phi})]$ for $i=1,\ldots
,31\,000$). Each interval matrix $[DF^{\phi}(\mathbf{I}_{i})]$ is computed
using (\ref{eq:DF-phi-I}). Thus we obtain a bound for $[DF^{\phi}(D_{\phi})]$
(displayed below with 3-digit rough accuracy; rounded up to ensure true
enclosure)\bigskip

\nopagebreak{ \noindent$\lbrack DF^{\phi}(D_{\phi})]=$ {\small
\begin{equation}
=\left(
\begin{array}
[c]{llll}%
\lbrack\mathtt{-0.0336,0.0335}] & [\mathtt{2.06,2.11}] &
[\mathtt{-0.0526,0.0521}] & [\mathtt{-0.0521,0.0526}]\\
\lbrack\mathtt{-2.15,-2.03}] & [\mathtt{-0.0422,0.0422}] &
[\mathtt{-0.0826,0.0827}] & [\mathtt{-0.0825,0.0827}]\\
\lbrack\mathtt{-0.0783,0.0782}] & [\mathtt{-0.0559,0.0566}] &
[\mathtt{2.43,2.64}] & [\mathtt{-0.0974,0.0962}]\\
\lbrack\mathtt{-0.0782,0.0783}] & [\mathtt{-0.0559,0.0566}] &
[\mathtt{-0.0962,0.0974}] & [\mathtt{-2.64,-2.43}]
\end{array}
\right)  \label{eq:DF-global-bound}%
\end{equation}
}}

We take $\alpha_{h}=\beta_{v}=2$ and $\alpha_{v}=\beta_{h}=\gamma=1,$ which
clearly satisfy (\ref{eq:alpha-beta-setup}). In our application we deal with a
single set $N_{p}=N_{0}=D_{\phi}$, which means that for this set $\rho=R.$
With our choice of parameters condition (\ref{eq:rho-cond-1}) clearly holds.

Based on (\ref{eq:DF-global-bound}) using (\ref{eq:est-expansion}%
-\ref{eq:est-epsilons}), (\ref{eq:kappa-h}) and (\ref{eq:kappa-v}) we compute
the constants $\kappa_{c}^{\text{forw}},$ $\kappa_{u}^{\text{forw}},$
$\kappa_{s}^{\text{forw}},$ $\kappa_{c}^{\text{back}},$ $\kappa_{u}%
^{\text{back}},$ $\kappa_{s}^{\text{back}},$ $\varepsilon_{u},$ $\varepsilon
_{s},$ $\delta^{u},$ $\delta^{s}$ needed for the verification of assumptions
of Theorem \ref{th:main}. The computed constants are written out in (\ref{eq:fin-num}) and (\ref{eq:kappa-bounds}) in the Appendix.

Finally, using the boxes $\mathbf{I}_{c,i}$ we also compute $[\pi_{x,y}%
F(\bar{B}_{c}^{R}\times\{0\}\times\{0\})]$ as the interval hull of all
$\left[  \pi_{x,y}F^{\phi}(\mathbf{I}_{c,i}\times\{0\}\times\{0\})\right]  $
for $i=1,\ldots,31\,000$ (displayed below with rough accuracy)%
\begin{align}
&  [\pi_{x,y}F(\bar{B}_{c}^{R}\times\{0\}\times\{0\})]\label{eq:F(BR,0)-bound}%
\\
&  =[\mathtt{-0.000954908,0.000819660}]\times\lbrack
\mathtt{-0.0009549120,0.000819658}],\nonumber
\end{align}
from which $E_{u},$ $E_{s}$ are computed using (\ref{eq:mth-as-cover-3}) (see (\ref{eq:fin-num}) in the Appendix). For
computation of each $\left[  \pi_{x,y}F^{\phi}(\mathbf{I}_{c,i}\times
\{0\}\times\{0\})\right]  $ we in fact need to further subdivide each box
$\mathbf{I}_{c,i}$ into nine parts (this is because $E_{u}$ and $E_{s}$ turn
out to be our most sensitive estimates). Based on all the computed constants
we verify assumptions (\ref{eq:mth-as-cone1}-\ref{eq:mth-as-cover-2}) of
Theorem \ref{th:main}.

The computer assisted part of the proof has taken 3 hours and 49 minutes of
computation on a standard laptop (it is possible to conduct much shorter proofs, but for less accurate enclosures of the manifold than above). Looking at the constants (\ref{eq:fin-num}), (\ref{eq:kappa-bounds}) written out in the
Appendix it is apparent that assumptions (\ref{eq:mth-as-cone1}),
(\ref{eq:mth-as-cone2}) of Theorem \ref{th:main} hold by a large margin. The
bottleneck lies in conditions (\ref{eq:mth-as-cover-1}) and
(\ref{eq:mth-as-cover-2}). This follows from the fact that the bounds computed
in (\ref{eq:F(BR,0)-bound}) are large in comparison to $r$ (see
(\ref{eq:mth-as-cover-3}) which binds the two together). This is because far away from the origin the 4-th order normal form no longer gives
an accurate enough estimate on the position of the manifold, and hence the
vector-field in the expansion/contraction direction becomes noticeably
nonzero. A simple remedy would be to use a higher order normal form, which would allow for obtaining a tighter enclosure and also a larger domain. This
would require longer computations and use of more capable hardware than a
standard laptop. Such computations though can easily be performed on clusters.

\begin{figure}[ptb]
\begin{center}
\includegraphics[width=12cm]{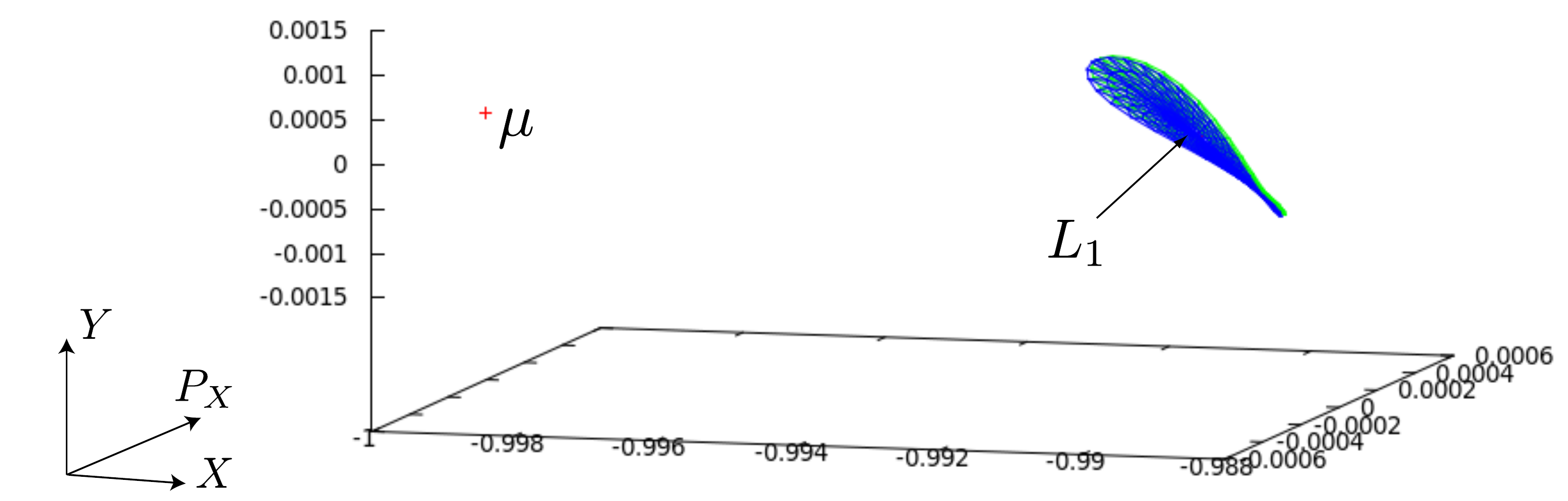}
\end{center}
\caption{A rough sketch of $\pi_{X,Y,P_{X}}\phi^{-1}(\bar{B}_{c}^{R-v}%
\times\bar{B}_{u}^{r}\times\bar{B}_{s}^{r})$, which gives us an idea of the
actual size and thickness of the investigated region in which we have proved
existence of the center manifold.}%
\label{fig:3dim}%
\end{figure}

\begin{figure}[ptb]
\begin{center}
\includegraphics[width=12cm]{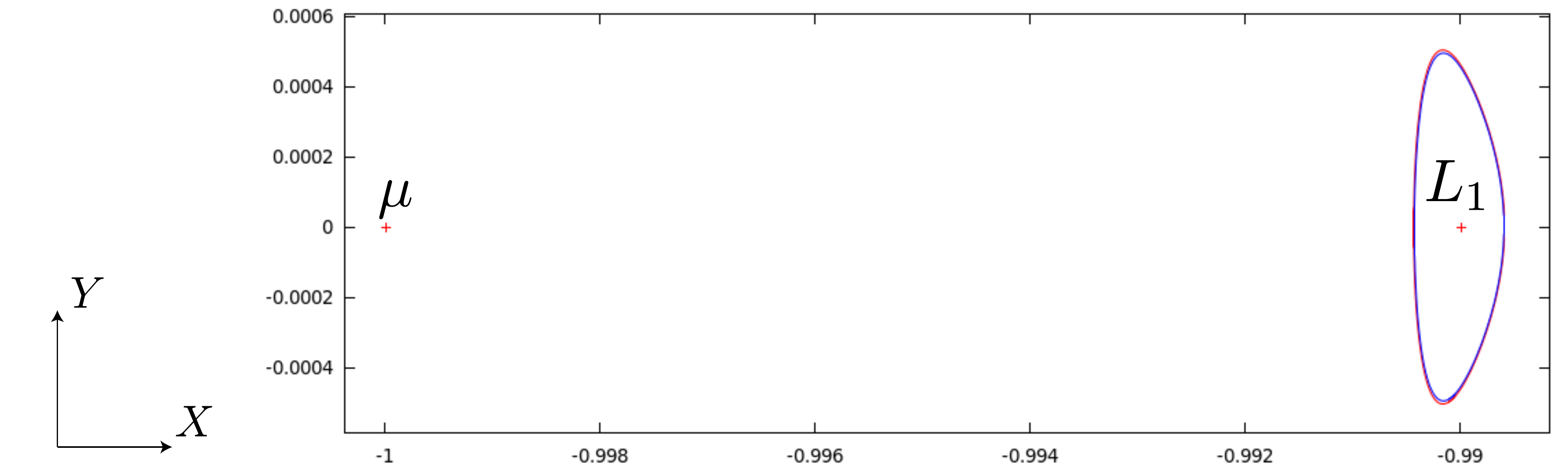}
\end{center}
\caption{A rough sketch of $\pi_{X,Y}\left(  \phi^{-1}(\partial \bar B_{c}^{R-v}
\times\bar B_{u}^{r} \times\bar B_{s}^{r}) \right)  $ in blue and
$\pi_{X,Y}\left(  \phi^{-1}(\partial \bar B_{c}^{R} \times\bar B_{u}^{r} \times\bar
B_{s}^{r}) \right)  $ in red. We have proved that the manifold is
contained in $\phi^{-1}(\bar B_{c}^{R-v} \times\bar B_{u}^{r} \times\bar
B_{s}^{r})$ and that orbits starting from it never leave $\phi^{-1}(\bar
B_{c}^{R} \times\bar B_{u}^{r} \times\bar B_{s}^{r})$ when going forwards or
backwards in time. }%
\label{fig:2dim}%
\end{figure}

Finally let us note that the size of the region in which the manifold is found
is not negligible. In Figures \ref{fig:3dim} and \ref{fig:2dim} we see our
region together with the smaller mass (Earth) in the original coordinates of
the system. Our set $D_{\phi}$ is a four dimensional "flattened disc", in Figure
\ref{fig:3dim} we can see that the disc is not too thick. On our plot the set
$\pi_{X,Y,P_{X}} \left(  \phi^{-1}\left(  D_{\phi}\right)  \right)  $ lies
between the two coloured flat discs (blue disc below, and green disc above; in this resolution they practically merge with one another).

\section{Closing remarks, future work}

\label{sec:concluding-remarks}

In this paper we have given a method for detection and proof of existence
of center manifolds in a practical domain of the system. We have successfully
applied the method to the Restricted Three Body Problem. The method is quite general. It can be applied to any system
with an integral of motion which allows for a computation of a normal form
around a fixed point. The method also works for arbitrary dimension, which
makes it a tool which can be applied to a large family of systems.

The strength of our approach lies in the fact that we can investigate and
prove existence of manifolds within large domains, and not only locally
around a fixed point. The weakness so far is that the method only establishes
Lipschitz type continuity of the manifold. In forthcoming work we plan to
remedy this deficiency. In our view, since we already have established Lipschitz continuity, similar tools combined with standard cohomology equation arguments can be applied to prove
higher order smoothness.

We would also like to mention that the method allows for rigorous enclosure of
the associated stable and unstable manifolds through cone conditions used in
the proof. This means that the presented method can be used as a starting
point for computation of foliations of stable/unstable manifolds, and next a scattering
map associated with splitting of separatrices. In our future work we plan to
conduct rigorous-computer-assisted computations of the scattering map for the
RTBP in the spirit of \cite{DelshamsMR08}. Such computations can then be used
in the study of structural stability or diffusion.

\section{Acknowledgements} We would like to express our thanks to Daniel Wilczak for frequent discussions and his assistance in the implementation of higher order computations in the CAPD library (http://capd.ii.uj.edu.pl).

\section{Appendix}

Here we list the bounds needed for the verification of assumptions of Theorem
\ref{th:main}. Below constants were computed using (\ref{eq:F(BR,0)-bound}),
(\ref{eq:DF-global-bound}) combined with (\ref{eq:mth-as-cover-3}),
(\ref{eq:est-expansion}), (\ref{eq:est-contraction}) and
(\ref{eq:est-epsilons})
\begin{equation}%
\begin{array}
[c]{lll}%
E_{u}=\mathtt{1.909815022732472,} &  & E_{s}=\mathtt{1.909823931307315,}\\
\delta_{u}=\mathtt{2.434904529896616,} &  & \delta_{s}%
=\mathtt{2.434911565550947,}\\
\varepsilon_{c}=\mathtt{0.09796031906285504,} &  & \varepsilon_{m}%
=\mathtt{0.09656707906887786,}\\
\varepsilon_{u}=\mathtt{0.09737656524499766,} &  & \varepsilon_{s}%
=\mathtt{0.09735689577043023.}%
\end{array}
\label{eq:fin-num}%
\end{equation}
Note that {$[DF^{\phi}(D_{\phi})]$ and }$[\pi_{x,y}F(\bar{B}_{c}^{R}%
\times\{0\}\times\{0\})]$ {in} (\ref{eq:DF-global-bound}),
(\ref{eq:F(BR,0)-bound}) are displayed with very rough accuracy. Above numbers
follow from their precise version from the CAPD software.

From (\ref{eq:dF-bound}) we have obtained the bounds $c^{u},c^{s}$ (see
(\ref{eq:est-central})) using the following simple estimates. Our matrix $\mathbf{C}$
from (\ref{eq:dF-bound}) is of the form (see (\ref{eq:DF-global-bound}))%
\[
\mathbf{C}=\left(
\begin{array}
[c]{cc}%
\boldsymbol{\varepsilon}_{1} & \mathbf{r}_{1}\\
\mathbf{r}_{2} & \boldsymbol{\varepsilon}_{2}%
\end{array}
\right)  .
\]
For any matrix $C=\left(
\begin{array}
[c]{cc}%
\varepsilon_{1} & r_{1}\\
r_{2} & \varepsilon_{2}%
\end{array}
\right)  \in\mathbf{C}$ and any $\theta=\left(  \theta_{1},\theta_{2}\right)
$ for which $\left\Vert \theta\right\Vert =1,$ using
\[
-\frac{1}{2}=-\frac{\theta_{1}^{2}+\theta_{2}^{2}}{2}\leq\theta_{1}\theta
_{2}\leq\frac{\theta_{1}^{2}+\theta_{2}^{2}}{2}=\frac{1}{2}%
\]
we have%
\begin{align}
&  \theta^{T}C\theta\nonumber\\
&  =\left(  r_{1}+r_{2}\right)  \theta_{2}\theta_{1}+\varepsilon_{1}\theta
_{1}^{2}+\varepsilon_{2}\theta_{2}^{2}\nonumber\\
&  \in\left[  -\max_{r_{1}\in\mathbf{r}_{1},r_{2}\in\mathbf{r}_{2}}%
\frac{|r_{1}+r_{2}|}{2}+\min_{\varepsilon_{i}\in\boldsymbol{\varepsilon}%
_{i},i=1,2}\varepsilon_{i},\max_{r_{1}\in\mathbf{r}_{1},r_{2}\in\mathbf{r}%
_{2}}\frac{|r_{1}+r_{2}|}{2}+\max_{\varepsilon_{i}\in\boldsymbol{\varepsilon
}_{i},i=1,2}\varepsilon_{i}\right]  .\label{eq:interval-c-bound}%
\end{align}
The bound (\ref{eq:interval-c-bound}) is easily computable using interval
arithmetic and (\ref{eq:DF-global-bound})
\begin{equation}%
\begin{array}
[c]{lll}%
c^{u}=\mathtt{0.08050236551044671,} &  & c^{s}=\mathtt{-0.08046115109310353.}%
\end{array}
\label{eq:fin-num-c}%
\end{equation}
Here once again the very rough rounding in (\ref{eq:DF-global-bound}) is
evident when compared with (\ref{eq:fin-num-c}).

Estimates (\ref{eq:fin-num}), (\ref{eq:fin-num-c}) give us%
\begin{equation}
\begin{array}
[c]{lll}%
\kappa_{c}^{\text{forw}}=\mathtt{\phantom{-}0.3233133031766185,} &  &
\kappa_{c}^{\text{back}}=\mathtt{-0.3232720887592754,}\\
\kappa_{u}^{\text{forw}}=\mathtt{\phantom{-}2.289103404031357,} &  &
\kappa_{u}^{\text{back}}=\mathtt{\phantom{-}2.191595652437820,}\\
\kappa_{s}^{\text{forw}}=\mathtt{-2.191592853354867,} &  & \kappa
_{s}^{\text{back}}=\mathtt{-2.289115357054329.}%
\end{array}\label{eq:kappa-bounds}
\end{equation}

\bibliographystyle{plain}



\end{document}